\documentclass[12pt]{article} 
\usepackage{amssymb} 
\usepackage{amscd}
\usepackage{amsthm}
\usepackage[utf8]{inputenc}
\usepackage{hyperref}
\newtheorem{theorem}{Theorem}

\newtheorem{proposition}[theorem]{Proposition}
\newtheorem{lemma}[theorem]{Lemma}
\newtheorem{corollary}[theorem]{Corollary}

\newtheorem{definition}{Definition}

\def\dim{{\mbox{dim}}}

\def\GL{{\mbox{GL}}}
\def\Hom {{\mbox{Hom}}}
\def\Ext {{\mbox{Ext}}}

\def\cali{{\cal I}}
 
\def\cala{{\cal A}} 
 
\def\calk{{\cal K}} 
\def\call{{\cal L}}
\def\calm{{\cal M}}

\def\calh{{\cal H}}
\def\cals{{\cal S}}

\def\calf{{\cal F}}

\def\calw{{\cal W}}

\def\bbbone{\mbox{\rm 1\hspace {-.6em} l}}

\def\ug{{\mathbf u}}

\def\Lg{{\mathbf L}}
\def\Rg{{\mathbf R}}

\def\Tor{{\mbox{Tor}}}

\begin{document}
\enlargethispage{3cm}

 \thispagestyle{empty}
\begin{center}
{\large\bf MULTILINEAR FORMS}
\end{center} 
\begin{center}
{\large\bf AND GRADED ALGEBRAS}
\end{center}
   
\vspace{0.3cm}

\begin{center} Michel DUBOIS-VIOLETTE
\footnote{Laboratoire de Physique Th\'eorique, UMR 8627, Universit\'e Paris XI,
B\^atiment 210, F-91 405 Orsay Cedex, France, 
Michel.Dubois-Violette$@$th.u-psud.fr\\
}
\end{center} \vspace{0,5cm}

\begin{abstract}
In this paper we investigate the class of the connected graded algebras which are finitely generated in degree 1, which are finitely presented with relations of degrees greater or equal to 2 and which are of finite global dimension $D$ and Gorenstein. For $D$ greater or equal to 4 we add the condition that these algebras are homogeneous and Koszul. It is shown that each such algebra is completely characterized by a multilinear form satisfying a twisted cyclicity condition and some other nondegeneracy conditions depending on the global dimension $D$. This multilinear form plays the role of a volume form and canonically identifies in the quadratic case to a nontrivial Hochschild cycle of maximal degree. Several examples including the Yang-Mills algebra and the extended 4-dimensional Sklyanin algebra are analyzed in this context. Actions of quantum groups are also investigated.

\end{abstract}

\newpage
\tableofcontents

\section{Introduction}

An important task which is at the beginning of noncommutative algebraic geometry is to provide good descriptions of the connected graded algebras which are finitely generated in degree 1, which are finitely presented with homogeneous relations of degrees $\geq 2$, which are of finite global dimension and which are Gorenstein (a generalization of the Poincar\'e duality property). These algebras play the role of homogeneous coordinates rings for the noncommutative versions of the projective spaces and, more generally, of algebraic varieties. It is usual to add the property of polynomial growth for these graded algebras \cite{art-sch:1987} but we refrain here to impose this property since it eliminates various interesting examples and plays no role in our arguments.\\
Some remarks are in order concerning the above class of algebras. The class of connected graded algebras which are finitely generated in degree 1 and finitely presented in degrees $\geq 2$ is a natural one for a generalization of the polynomial algebras. Concerning the global dimension it is an important fact which is well known \cite{art-tat-vdb:1990} that for this class of algebras it coincides with the projective dimension of the trivial module (the ground field) and it has been shown recently  \cite{ber:2005}
that it also coincides with the Hochschild dimension. Thus for these algebras it is {\sl the dimension} from the homological point of view and the requirement of finite dimensionality is clear. That the Gorenstein property is  a generalization of the Poincar\'e duality property is already visible if one thinks of the minimal projective resolution of the trivial module as an analog of differential forms and this has been made precise at the Hochschild homological level in \cite{ber-mar:2006} (see also \cite{vdb:1998},\cite{vdb:2002}) . \\
In this paper we shall restrict attention to the smaller class of the algebras which are also homogeneous and Koszul. Homogeneous means that all the relations are of the same degree, say $N\geq 2$ and we speak then of homogeneous algebras of degree $N$ or of $N$-homogeneous algebras. For homogeneous algebras the notion of Koszulity has been introduced in \cite{ber:2001a} and various notions such as the Koszul duality, etc. generalizing the ones occuring in the quadratic case \cite{pri:1970}, \cite{man:1988} have been introduced in \cite{ber-mdv-wam:2003}. It should be stressed that the Koszul property is really a desired property \cite{man:1988} and this is the very reason why we restrict attention to homogeneous algebras since it is only for these algebras that we know how to formulate this property for the moment. It is worth noticing here that these restrictions are immaterial in the case of global dimension $D=2$ and $D=3$. Indeed,  as pointed out in \cite{ber-mar:2006} any connected graded algebra which is finitely generated in degree 1, finitely presented with relations of degree $\geq 2$ and which is of global dimension $D=2$ or $D=3$ and  Gorenstein, is $N$-homogeneous and Koszul with $N=2$ for $D=2$ and $N\geq 2$ for $D=3$. However this is already not the case for the global dimension $D=4$ \cite{art-sch:1987}, \cite{lu-pal-wu-zha:2007}.\\
In the following we shall give detailed proofs of  results announced in 
\cite{mdv:2005} which allow to identify the moduli space of the algebras $\cala$ as above with the moduli space of multilinear forms $w$ with specific properties. For each $\cala$, the multilinear form $w$ or more precisely $\bbbone\otimes w$ ($\bbbone \in \cala$) plays the role of a volume element in the Koszul resolution of the trivial $\cala$-module. It turns out that this is also true from the point of view of the Hochschild homology of $\cala$ at least in the quadratic case as will be shown. This gives another bridge, besides the deep ones described in \cite{ac-mdv:2003} and 
\cite{ac-mdv:2008}, between noncommutative differential geometry (\cite{ac:1986a},\cite{ac:1994})  and noncommutative algebraic geometry (\cite{staf:2002} and references therein). We shall analyse several examples in order to illustrate the concepts introduced throughout the paper. Finally we shall introduce related Hopf algebras.\\

It is worth noticing here that in \cite{bon-pol:1994} it has been already shown that the quadratic algebras which are Koszul of finite global dimension and Gorenstein are determined by multilinear forms. This is of course directly related with the results of the present paper and we shall come back to this point later (Section 9).\\

The plan of the paper is the following. In Section 2, we investigate the case of the global dimension $D=2$ and we show that the algebras of the relevant class are associated with the nondegenerate bilinear forms and correspond to the natural quantum spaces for the action of the quantum groups of the associated nondegenerate bilinear forms \cite{mdv-lau:1990}. In Section 3, we introduce and discuss the concept of preregular multilinear form. It turns out that, as shown in this paper (in Section 5), all homogeneous Koszul-Gorenstein algebras of finite global dimension are associated with preregular multilinear forms satisfying some other regularity conditions depending on the global dimension $D$. The case $D=3$ is analysed in Section 4. In Section 5, we define and study the homogeneous algebras associated with multilinear forms. Section 6 consists of the analysis of several examples which illustrate the different items of the paper. In Section 7 the semi-cross product is investigated for the above class of algebras (introduced in Section 5). In Section 8 we define quantum groups preserving the multilinear forms which act on the quantum spaces corresponding to the homogeneous algebras associated with these multilinear forms generalizing thereby the situation for $D=2$ described in Section 2. In Section 9 we discuss several important points connected with the present formulation. For the reader convenience we have added an appendix on homogeneous algebras and an appendix on the quantum group of a nondegenerate bilinear form at the end of the paper.\\
Throughout the paper $\mathbb K$ denotes a field which we assume to be algebraically closed of characteristic zero (though most of our results are independent of this assumption) and all algebras and vector spaces are over $\mathbb K$. The symbol $\otimes$ denotes the tensor product over $\mathbb K$ and if $E$ is a $\mathbb K$-vector space, $E^\ast$ denotes its dual vector space. We use the Einstein summation convention of repeated up down indices in the formulas.

\section{Bilinear forms and global dimension $D=2$}

Let $b$ be a nondegenerate bilinear form on $\mathbb K^{s+1}$ ($s\geq 1$) with components $B_{\mu\nu}=b(e_\mu,e_\nu)$ in the canonical basis $(e_\lambda)_{\lambda\in \{0,\dots,s\}}$ of $\mathbb K^{s+1}$ and let $\cala=\cala(b,2)$ be the quadratic algebra generated by the elements $x^\lambda$ ($\lambda\in \{0,\dots,s\}$) with the relation
\begin{equation}
B_{\mu\nu}x^\mu x^\nu=0
\end{equation}
that is, using the notations of \cite{ber-mdv-wam:2003} (see Appendix 1), one has $\cala=A(E,R)$ with $E=\oplus_\lambda\mathbb K x^\lambda=\cala_1$ and $R=\mathbb K\ B_{\mu\nu}x^\mu\otimes x^\nu\subset E^{\otimes^2}$.
\begin{lemma}\label{reg}
Let $a$ be an element of degree 1 of $\cala$ with $a\not= 0$ $($i.e. $a\in E\backslash \{0\})$. Then $ya=0$ or $ay=0$ for $y\in \cala$ implies $y=0$.
\end{lemma}
\noindent \underbar{Proof}. One has
\begin{equation}
a=a_\lambda x^\lambda
\end{equation}
with $(a_\lambda)\not= 0$ in $\mathbb K^{s+1}$. On the other hand, by the very definition of $\cala$ by generators and relation, $ya=0$ (resp. $ay=0$) is equivalent to
\begin{equation}
ya_\lambda=zx^\mu\  B_{\mu\lambda}\ \ \ (\mbox{resp.} ya_\lambda=x^\mu z\  B_{\lambda\mu})
\label{eq2.3}
\end{equation}
$(\lambda\in \{0,\dots,s\})$ for some $z\in \cala$. Let $y$ be in $\oplus^n_{k=0}\cala_k$, we shall prove the statement by induction on $n$. For $n=0$ (i.e. $y$ of degree 0) the statement is clear ($\cala$ is connected). Assume that the statement is true for $n\leq  p$ and let $y\in \oplus^{p+1}_{k=0}\cala_k$ be such that $ya=0$ (resp. $ay=0$). Then $z\in \oplus^p_{k=0}\cala_k$ in (\ref{eq2.3}) and one has $zx^\mu B_{\mu\lambda} v^\lambda=0$  (resp. $B_{\lambda\mu} v^\lambda x^\mu z=0)$ for $(v^\lambda)\not=0$ in $\mathbb K^{s+1}$ such that $a_\lambda v^\lambda=0$ in $\mathbb K$. So, replacing $a$ by $x^\mu B_{\mu\lambda} v^\lambda$ (resp. $B_{\lambda \mu} v^\lambda x^\mu$), one has $z=0$ by the induction hypothesis. This implies $y=0$ by (\ref{eq2.3}) since $(a_\lambda)\not=0$. $\square$

The matrix $B=(B_{\mu\nu})$ of the components of $b$ is invertible and we denote by $B^{\mu\nu}$ the matrix elements of its inverse that is one has 
\begin{equation}
B^{\lambda\rho}B_{\rho\mu}=\delta^\lambda_\mu
\end{equation}
$\lambda,\mu\in \{0,\dots,s\}$. The $B^{\mu\nu}$ are the component of a bilinear form on the dual of $\mathbb K^{s+1}$ in the dual basis of the canonical basis $(e_\lambda)$. Notice that with the definitions above the vector space $E$ identifies canonically with the dual of $\mathbb K^{s+1}$ while the $x^\lambda$ identify with the elements of the dual basis of the canonical basis $(e_\lambda)$ of $\mathbb K^{s+1}$. These identifications allow for instance to write $b\in E^{\otimes^2}$ since the involved vector spaces are finite-dimensional.\\
Let us investigate the structure of the dual quadratic algebra $\cala^!$ of $\cala$ \cite{man:1988}, \cite{ber-mdv-wam:2003}. Letting $E^\ast=\mathbb K^{s+1}$ be the dual vector space of $E$, one has $\cala^!=A(E^\ast,R^\perp)$ where $R^\perp \subset (E^\ast)^{\otimes^2}=(E^{\otimes^2})^\ast$ is the orthogonal of $R=\mathbb K B_{\mu\nu}x^\mu\otimes x^\nu\subset E^{\otimes^2}$.
\begin{lemma}\label{d2}
The dual quadratic algebra of $\cala$ is the quadratic algebra $\cala^!$ generated by the elements $e_\lambda$ $(\lambda\in \{0,\dots,s\})$ with the relations
\begin{equation}
e_\mu e_\nu=\frac{1}{s+1}B_{\mu\nu} B^{\tau\rho}e_\rho e_\tau
\end{equation}
for $\mu,\nu\in \{0,\dots,s\}$. One has $\cala^!_0=\mathbb K\bbbone \simeq \mathbb K$, $\cala^!_1=E^\ast=\oplus_\lambda\mathbb K e_\lambda\simeq \mathbb K^{s+1}$, $\cala^!_2=\mathbb K B^{\mu\nu}e_\nu e_\mu\simeq \mathbb K$ and $\cala^!_n=0$ for $n\geq 3$.
\end{lemma}

\noindent\underbar{Proof}. One has $\langle e_\mu\otimes e_\nu-\frac{1}{s+1} B_{\mu\nu} B^{\tau \rho}e_\rho \otimes e_\tau,\ \ B_{\lambda\sigma}x^\lambda\otimes x^\sigma\rangle=0$ so the $e_\mu \otimes e_\nu-\frac{1}{s+1}B_{\mu\nu}B^{\tau\rho}e_\rho \otimes e_\tau$ are in $R^\perp$ and it is not difficult to see that they span $R^\perp$ which proves the first part of the lemma including the identifications of $\cala^!_0,\cala^!_1$ and $\cala^!_2$. It remains to show that $\cala^!_3=0$. Setting $\xi=B^{\alpha\beta}e_\beta e_\alpha$ for the generator of $\cala^!_2$ one has $(e_\lambda e_\mu)e_\nu=\frac{1}{s+1}B_{\lambda\mu}\xi e_\nu$ which (by contraction with $B^{\nu\mu}$) implies $e_\lambda\xi=\frac{1}{s+1} B^{\nu\mu} B_{\lambda\mu}\xi e_\nu$ while one also has $e_\lambda(e_\mu e_\nu)=e_\lambda\frac{1}{s+1}B_{\mu\nu}\xi$ which (by contraction with $B^{\mu\lambda}$) implies $\xi e_\nu=\frac{1}{s+1} B^{\mu\lambda}B_{\mu\nu} e_\lambda\xi$. It follows that one has $e_\lambda\xi=\left(\frac{1}{s+1}\right)^2 e_\lambda\xi$ that is $e_\lambda\xi=0$ since $s\geq 1$ and thus $e_\lambda e_\mu e_\nu=0$ for $\lambda,\mu,\nu \in \{0,\dots,s\}$ which means $\cala^!_3=0$. $\square$

In view of this lemma the Koszul complex of $\cala$
\[
\dots \rightarrow \cala \otimes \cala^{!\ast}_{n+1} \stackrel{d}{\rightarrow} \cala \otimes \cala^{!\ast}_n\rightarrow \dots
\]
reads here
\begin{equation}
0\rightarrow \cala \stackrel{x^tB}{\rightarrow} \cala^{s+1} \stackrel{x}{\rightarrow} \cala \rightarrow 0
\end{equation}
where $\cala^{s+1}=(\cala,\dots,\cala)$, $x$ means right multiplication by the column $(x^\lambda)$ and $x^tB$ means right multiplication by the row $(x^\mu B_{\mu\lambda})$. Lemma \ref{reg} implies that $\cala \stackrel{x^tB}{\rightarrow} \cala^{s+1}$ is injective while the definition of $\cala$ by generators and relation means that the sequence $\cala\stackrel{x^tB}{\rightarrow}\cala^{s+1}\stackrel{x}{\rightarrow}\cala\stackrel{\varepsilon}{\rightarrow} \mathbb K \rightarrow 0$ is exact ($\varepsilon$ being the projection on degree 0).  Therefore $\cala$ is Koszul of global dimension 2 and the exact sequence of left $\cala$-module
\begin{equation}
0\rightarrow \cala \stackrel{x^tB}{\rightarrow} \cala^{s+1} \stackrel{x}{\rightarrow}\cala \stackrel{\varepsilon}{\rightarrow} \mathbb K \rightarrow 0
\label{eq2.7}
\end{equation}
is the Koszul (free) resolution of the trivial left $\cala$-module $\mathbb K$. By transposition and by using the invertibility of $B$, it follows that $\cala$ is also Gorenstein.\\
Conversely let $\cala$ be a connected graded algebra generated by $s+1$ elements $x^\lambda$ of degree 1 with (a finite number of) relations of degrees $\geq 2$ which is of global dimension 2 and  Gorenstein. Then, as pointed out in the introduction, it is known (and easy to show) that $\cala$ is quadratic and Koszul. The Gorenstein property implies that the space of relations $R$ is 1-dimensional so $\cala$ is generated by the $x^\lambda$ with relation $B_{\mu\nu}x^\mu x^\nu=0$ ($B_{\mu\nu}\in \mathbb K$, $\mu,\nu\in \{0,\dots,s\}$) and the Koszul resolution of $\mathbb K$ is of the above form (\ref{eq2.7}). Furthermore the Gorenstein property also implies that $B$ is invertible so the corresponding bilinear form $b$ on $\mathbb K^{s+1}$ is nondegenerate and $\cala$ is of the above type (i.e. $\cala=\cala(b,2)$). This is summarized by the following theorem.
\begin{theorem}\label{KGD2}
Let $b$ be a nondegenerate bilinear form on $\mathbb K^{s+1}$ $(s\geq 1)$ with components $B_{\mu\nu}=b(e_\mu,e_\nu)$ in the canonical basis $(e_\lambda)$ of $\mathbb K^{s+1}$. Then the quadratic algebra $\cala$ generated by the elements $x^\lambda$ $(\lambda\in \{ 0,\dots,s\})$ with the relation $B_{\mu\nu}x^\mu x^\nu=0$ is Koszul of global dimension 2 and Gorenstein. Furthermore any connected graded algebra generated by $s+1$ element $x^\lambda$ of degree 1 with relations of degree $\geq 2$ which is of global dimension 2 and Gorenstein is of the above kind for some nondegenerate bilinear form $b$ on $\mathbb K^{s+1}$. 
\end{theorem}

There is a canonical right action $b\mapsto b\circ L$ ($L\in GL(s+1,\mathbb K)$) of the linear group on bilinear forms, where
\begin{equation}
(b\circ L)(X,Y)=b(LX,LY)
\end{equation}
for $X,Y\in \mathbb K^{s+1}$, which preserves the set of nondegenerate bilinear forms and one has the following straightforward result which is worth noticing in comparison with the similar one in global dimension $D=3$ which is less obvious (see in Section 4).
\begin{proposition} \label{M2}
Two nondegenerate bilinear forms $b$ and $b'$ on $\mathbb K^{s+1}$ correspond to isomorphic graded algebras $\cala(b,2)$ and $\cala(b',2)$ if and only if they belong to the same $GL(s+1,\mathbb K)$-orbit, i.e. if $b'=b\circ L$ for some $L\in GL(s+1,\mathbb K).$\end{proposition}
In view of Theorem \ref{KGD2} and Proposition \ref{M2} it is natural to define the {\sl moduli space} $\calm_s(2)$ of the quadratic algebras with $s+1$ generators which are Koszul of global dimension 2 and Gorenstein to be the space of $GL(s+1,\mathbb K)$-orbits of nondegenerate bilinear forms on $\mathbb K^{s+1}$. The {\sl moduli space} $\calm(2)$ of the connected graded algebras which are finitely generated in degree 1 and finitely presented with relations of degrees $\geq 2$ and which are of global dimension 2 and Gorenstein being then the (disjoint) union $\calm(2)=\cup_{s\geq 1} \calm_s(2)$.\\

The Poincar\'e series $P_\cala(t)=\sum_n \dim(\cala_n)t^n$ of a graded algebra $\cala$ as above is given by \cite{pri:1970}, \cite{man:1988}
\begin{equation}
P_\cala(t)=\frac{1}{1-(s+1)t+t^2}
\end{equation}
which implies exponential growth for $s\geq 2$. For $s=1\ \ (s+1=2)$ the algebra has polynomial growth so it is regular in the sense of \cite{art-sch:1987}. In the latter case it is easy to classify the $GL(2,\mathbb K)$-orbits of nondegenerate bilinear forms on $\mathbb K^2$ according to the rank $\mathbf{rk}$ of their symmetric part \cite{mdv-lau:1990}:
\begin{itemize}
\item
$\mathbf{rk}=0$ - there is only one orbit which is the orbit of $b=\varepsilon$ with matrix of components $B=\left(\begin{array}{cc}
0 & -1\\
1 & 0
\end{array}\right)$
which corresponds to the relation\linebreak[4] $x^1x^2-x^2x^1=0$ i.e. to the polynomials algebra $\mathbb K[x^1,x^2]$,

\item
$\mathbf{rk}=1$ - there is only one orbit which is the orbit of $b$ with matrix of components $B=\left(\begin{array}{cc}
0 & -1\\
1 & 1
\end{array}\right)$ which corresponds to the relation \linebreak[4] $x^1x^2-x^2x^1-(x^2)^2=0$, 
\item
$\mathbf{rk}=2$ - the orbits are the orbits of $b=\varepsilon_q$ with matrix of components $B=\left(\begin{array}{cc}
0 & -1\\
q & 0
\end{array}\right)$ for $q\not=0$ and $q\not=1$ ($q^2\not=q$) modulo $q\sim q^{-1}$ which corresponds to the relation $x^1x^2-qx^2x^1=0$.
\end{itemize}
 Thus for $s+1=2$ one recovers the usual description of regular algebras of global dimension 2  \cite{irv:1979}  \cite{art-sch:1987}.\\

The algebra $\cala_q$ of the latter case $\mathbf{rk}=2$ corresponds to the Manin plane which is the natural quantum space for the action of the quantum group $SL_q(2, \mathbb K)$ \cite{man:1988}. More generally, given $s\geq 1$ and a nondegenerate bilinear form $b$ on $\mathbb K^{s+1}$ (with matrix of components $B$), the algebra $\cala$ of Theorem \ref{KGD2} corresponds to the natural quantum space for the action of the quantum group of the nondegenerate bilinear form $b$ \cite{mdv-lau:1990} (see Appendix 2). The complete analysis of the category of representations of this quantum group has been done in \cite{bic:2003b}.\\

A lot of general results on regular rings of dimension 2 is contained in Reference \cite{zha:1998}. In particular, Theorem 3 above is covered par Theorem 0.1 of \cite{zha:1998} while Theorem 0.2 (2) of \cite{zha:1998} stating that $\cala$ is a domain clearly implies Lemma 1 above.\\
In Reference \cite{ber:2009}, the cases where the bilinear form $b$ is degenerate is investigated and, more generally the case of the graded algebras $\cala$ having a single homogeneous relation is analized there.

\section{Multilinear forms}

In this section $V$ is a vector space with $\dim(V)\geq 2$ and $m$ is an integer with $m\geq 2$.

\begin{definition}\label{tc}
Let $Q$ be an element of the linear group $GL(V)$. A $m$-linear form $w$ on $V$ will be said to satisfy the $Q$-twisted cyclicity condition or simply to be $Q$-cyclic if one has
\begin{equation}
w(X_1,\dots,X_m)=w(QX_m,X_1,\dots,X_{m-1})
\end{equation}
for any $X_1,\dots,X_m \in V$.
\end{definition}
Let $w$ be $Q$-cyclic then one has
\begin{equation}
w(X_1,\dots,X_m)=w(QX_k,\dots,QX_m,X_1,\dots,X_{k-1})
\end{equation}
for any $2\leq k\leq m$ and finally
\begin{equation}
w(X_1,\dots,X_m)=w(QX_1,\dots, QX_m)
\end{equation}
for any $X_1,\dots, X_m\in V$ which means that $w$ is invariant by $Q$, $w=w\circ Q$.\\
Let now $w$ be an arbitrary $Q$-invariant ($w=w\circ Q$) $m$-linear form on $V$, then $\pi_Q(w)$ defined by
\begin{equation}
m\pi_Q(w)(X_1,\dots, X_m)=w(X_1,\dots,X_m)+\sum^m_{k=2}w(QX_k,\dots, QX_m,X_1,\dots, X_{k-1})
\end{equation}
is a $Q$-cyclic $m$-linear form on $V$ and the linear mapping $\pi_Q$ is a projection of the space of $Q$-invariant $m$-linear forms onto the subspace of the $Q$-cyclic ones ($(\pi_Q)^2=\pi_Q$).\\
Notice that to admit a non zero $Q$-invariant multilinear form is a nontrivial condition on $Q\in GL(V)$ since it means that the operator $w\mapsto w\circ Q$ has an eigenvalue equal to 1. For instance there is no non zero $(-\bbbone)$-invariant $m$-linear form if $m$ is odd.\\
Let us consider the right action $w\mapsto w\circ L$ ($L\in GL(V)$) of the linear group on the space of $m$-linear forms on $V$. If $w$ is $Q$-invariant then $w\circ L$ is $L^{-1}QL$-invariant and if $w$ is $Q$-cyclic then $w\circ L$ is $L^{-1}QL$-cyclic.

\begin{definition} \label{PR}
A $m$-linear form $w$ on $V$ will be said to be preregular if it satisfies the conditions $(i)$ and $(ii)$ below.\\
$(i)$ If $X\in V$ satisfies $w(X,X_1,\dots,X_{m-1})=0$ for any $X_1,\dots,X_{m-1}\in V$, then $X=0$.\\
$(ii)$ There is a $Q_w\in GL(V)$ such that $w$ is $Q_w$-cyclic.
\end{definition}
The condition $(i)$ implies that the element $Q_w$ of $GL(V)$ such that $(ii)$ is satisfied is unique for a preregular  $m$-linear form $w$ on $V$. In view of $(ii)$ a preregular $w$ is such that if $X\in V$ satisfies 
\[
w(X_1,\dots, X_k,X,X_{k+1},\dots,X_{m-1})=0
\]
 for any $X_1,\dots, X_{m-1}\in V$ then $X=0$. A $m$-linear form $w$ satisfying this latter condition for any $k$ $ (0\leq k\leq m)$ will be said to be 1-{\sl site nondegenerate}. Thus a preregular $m$-linear form is a 1-site nondegenerate twisted cyclic $m$-linear form.\\
The set of preregular $m$-linear forms on $V$ is invariant by the action of the linear group $GL(V)$ and one has
\begin{equation}
Q_{w\circ L}=L^{-1}Q_wL,\ \ \ \forall L\in GL(V)
\end{equation}
for a preregular $m$-linear form $w$ on $V$.\\
A bilinear form $b$ on $\mathbb K^{s+1}$ ($s\geq 1$) is preregular if and only if it is nondegenerate. Indeed if $b$ is preregular, it is nondegenerate in view of $(i)$. Conversely if $b$ is nondegenerate with matrix of components $B$ then one has $b(X,Y)=b(Q_bY,X)$ with $Q_b=(B^{-1})^tB$.\\
As pointed out in the introduction and as will be shown later all homogeneous Koszul algebras of finite global dimension which are also Gorenstein are associated with preregular multilinear forms satisfying some other conditions depending on the global dimension $D$. Let us spell out a condition for the case $D=3$ which will be the object of the next section.

\begin{definition}\label{3R}
Let $N$ be an integer with $N\geq 2$. A $(N+1)$-linear form $w$ on $V$ will be said to be 3-regular if it is preregular and satisfies the following condition $(iii)$.\\
$(iii)$ If $L_0$ and $L_1$ are endomorphisms of $V$ satisfying 
\[
w(L_0X_0,X_1,X_2,\dots,X_N)=w(X_0,L_1X_1,X_2,\dots,X_N)
\]
for any $X_0,\dots,X_N\in V$, then $L_0=L_1=\lambda\bbbone$ for some $\lambda\in \mathbb K$.
\end{definition}
The set of all 3-regular $(N+1)$-linear forms on $V$ is also invariant by the right action of $GL(V)$.\\
Notice that condition $(iii)$ is a sort of two-sites (consecutive, etc.) nondegenerate condition. Consider the following more natural two-sites condition $(iii)'$ for a $(N+1)$-linear form $w$ on $V(N\geq 2)$.\\

\noindent $(iii)'$ {\sl If} $\sum_iY_i\otimes Z_i \in V\otimes V$ {\sl satisfies}
\[
\sum_i w(Y_i,Z_i,X_1,\dots,X_{N-1})=0
\]
{\sl for any} $X_1,\dots,X_{N-1} \in V$, {\sl then} $\sum_iY_i\otimes Z_i=0$.\\

It is clear that the condition $(iii)'$ implies the condition $(iii)$, but it is a strictly stronger condition. For instance take $V=\mathbb K^{N+1}$ and let $\varepsilon$ be the completely antisymmetric $(N+1)$-linear form on $\mathbb K^{N+1}$ such that $\varepsilon(e_0,\dots e_N)=1$ in the canonical basis $(e_\alpha)$ of $\mathbb K^{N+1}$. Then $\varepsilon$ is 3-regular but does not satisfy $(iii)'$ since for $Y\otimes Z + Z \otimes Y \in V\otimes V$ one has
\[
\varepsilon (Y,Z,X_1,\dots,X_{N-1})+\varepsilon (Z,Y,X_1,\dots,X_{N-1})=0
\]
identically. 

\section{ Global dimension $D=3$}
Let $w$ be a preregular $(N+1)$-linear form on $\mathbb K^{s+1}$ (with $N\geq 2$ and $s\geq 1$) with components $W_{\lambda_0\dots \lambda_N}=w(e_{\lambda_0},\dots,e_{\lambda_N}$) in the canonical basis $(e_\lambda)$ of $\mathbb K^{s+1}$ and let $\cala=\cala(w,N)$ be the $N$-homogeneous algebra generated by the $s+1$ elements $x^\lambda$ ($\lambda\in \{0,\dots,s\}$) with the $s+1$ relations
\begin{equation}
W_{\lambda\lambda_1\dots \lambda_N}x^{\lambda_1} \dots x^{\lambda_N}=0\ \ \ (\lambda\in \{0,\dots,s\})\label{4.1}
\end{equation}
that is, again with the notations of \cite{ber-mdv-wam:2003},  one has $\cala=A(E,R)$ with $E=\oplus_\lambda\mathbb K x^\lambda=\cala_1$ and $R=\sum_\lambda\mathbb K W_{\lambda\lambda_1\dots \lambda_N}x^{\lambda_1}\otimes \dots \otimes x^{\lambda_N} \subset E^{\otimes^N}$. Notice that the condition $(i)$ implies that the latter sum is a direct sum i.e. $\dim(R)=s+1$. \\

\noindent \underbar{Remark}. Since $w$ is 1-site nondegenerate, there is a (non unique) $(N+1)$-linear form $\tilde w$ on the dual vector space of $\mathbb K^{s+1}$ (i.e. on $E$) with components $\tilde W^{\lambda_0\dots \lambda_N}$ in the dual basis of $(e_\lambda)$ such that $\tilde W^{\mu \lambda_1\dots\lambda_N}W_{\lambda_1\dots \lambda_N\nu}=\delta^\mu_\nu$. Let $\theta_\lambda$ ($\lambda\in \{0,\dots,s\}$) be the generators of $\cala^!=\cala(w,N)^!$ corresponding to the $x^\lambda$ (dual basis). Then the $\Theta^\lambda =\tilde W^{\lambda \lambda_1\dots \lambda_N}\theta_{\lambda_1}\dots \theta_{\lambda_N}$ ($\lambda\in \{0,\dots,s\}$) form a basis of $\cala^!_N$, the relations of $\cala^!$ read
\[
\theta_{\mu_1}\dots \theta_{\mu_N}=W_{\mu_1\dots \mu_N \lambda} \Theta^\lambda, \ \ \ (\mu_k\in \{0,\dots,s\}).
\]
and the $\Theta^\lambda$ are independent of the choice of $\tilde w$ as above. If furthermore $w$ is 3-regular then Proposition \ref{CK3} in Section 9 implies that $\cala^!_{N+1}$ is 1-dimensional generated by $\Theta^\lambda\theta_\lambda$ and that $\cala^!_n=0$ for $n\geq N+2$.\\
Notice that, in view of the $Q_w$-cyclicity, the relations (\ref{4.1}) of $\cala$ read as well $W_{\mu_1\dots \mu_N\lambda} x^{\mu_1}\dots x^{\mu_N}=0$, ($\lambda\in \{0,\dots,s\}$).

\begin{theorem}\label{G3} Let $\cala$ be a connected graded algebra which is finitely generated in degree 1, finitely presented with relations of degree $\geq 2$ and which is of global dimension $D=3$ and Gorenstein. Then $\cala=\cala(w,N)$ for some 3-regular $(N+1)$-linear form $w$ on $\mathbb K^{s+1}$.
\end{theorem}

\noindent \underbar{Proof}. As pointed out in \cite{ber-mar:2006} (see also in the introduction) $\cala$ is $N$-homogeneous with $N\geq 2$ and is Koszul. It follows then from the general theorem \ref{KGD} of next section that $\cala=\cala(w,N)$ for some preregular $(N+1)$-linear form on $\mathbb K^{s+1}$. Let us show that $w$ is in fact 3-regular.\\
The Koszul resolution of the trivial left $\cala$-module $\mathbb K$ reads
\begin{equation}
0\rightarrow \cala\otimes w \stackrel{d}{\rightarrow} \cala\otimes R \stackrel{d^{N-1}}{\rightarrow} \cala\otimes E \stackrel{d}{\rightarrow} \cala \rightarrow \mathbb K \rightarrow 0
\label{KR3}
\end{equation} where $E=\oplus_\mu \mathbb K x^\mu$, $R=\oplus_\mu \mathbb K W_{\mu \mu_1\dots \mu_N} x^{\mu_1}\otimes\dots \otimes x^{\mu_N}\subset E^{\otimes^N}$ and where $w$ is identified with the element $W_{\mu_0 \dots \mu_N}x^{\mu_0}\otimes \dots \otimes x^{\mu_N}$ of $E^{\otimes^{N+1}}$; the $N$-differential $d$ being induced by the mapping $a\otimes (x^0\otimes x^1\otimes \dots \otimes x^n)\mapsto ax^0\otimes (x^1\otimes \dots \otimes x^n)$ of $\cala\otimes E^{\otimes^{n+1}}$ into $\cala\otimes E^{\otimes^n}$, \cite{ber-mdv-wam:2003} (see Appendix 1).\\

Assume that the matrices $L_0,L_1\in M_{s+1}(\mathbb K)$ are such that one has
\begin{equation}
L_{0\ \mu_0}^\mu W_{\mu\mu_1\dots \mu_N}=L_{1\ \mu_1}^\mu W_{\mu_0\mu\mu_2\dots \mu_N} \label{4.2}
\end{equation}
and let $a^\mu\in \cala_1$ be the elements $a^\mu=L^\mu_{1\ \nu} x^\nu$. Equation (\ref{4.2}) implies $W_{\mu_0\mu\mu_2\dots\mu_N} a^\mu x^{\mu_2}\dots x^{\mu_N}=0$ or equivalently in view of Property $(ii)$ of Definition \ref{PR}
\[
W_{\mu\mu_1\dots \mu_N}a^\mu x^{\mu_1}\dots x^{\mu_{N-1}}\otimes x^{\mu_N}=0
\]
which also reads
\begin{equation}
d^{N-1}(a^\mu\otimes W_{\mu\mu_1\dots \mu_N}x^{\mu_1}\otimes \dots \otimes^{\mu_N})=0
\end{equation}
for the element $a^\mu\otimes W_{\mu\mu_1\dots \mu_N}x^{\mu_1}\otimes \dots \otimes x^{\mu_N}$ of $\cala_1\otimes R$. Exactness of the sequence \ref{KR3} at $\cala\otimes R$ implies that one has
\begin{equation}
a^\mu\otimes W_{\mu\mu_1\dots \mu_N}x^{\mu_1}\otimes \dots \otimes x^{\mu_N}=d(\lambda \bbbone\otimes w)
\end{equation}
for some $\lambda\in \mathbb K$. This implies
\[
a^\mu=L^\mu_{1\ \nu}x^\nu=\lambda x^\mu
\]
in view of Property $(i)$ of Definition \ref{PR}. Using again Property $(i)$, one finally obtains
\[
L_0=L_1=\lambda\bbbone \ \  (\in M_{s+1}(\mathbb K))
\]
as consequence of (\ref{4.2}) which means that $w$ is 3-regular.$\square$\\

The Poincar\'e series of such a $N$-homogeneous algebra $\cala$ which is Koszul of global dimension 3 and which is Gorenstein is given by \cite{art-tat-vdb:1991},\cite{mdv-pop:2002}
\begin{equation}
P_\cala(t)=\frac{1}{1-(s+1)t+(s+1)t^N-t^{N+1}}
\label{S3}
\end{equation}
where $s+1=\dim (\cala_1)$ is as before the number of independent generators (of degree 1). It follows from this formula that $\cala$ has exponential growth if $s+1+N>5$ while the case $s+1=2$ and $N=2$ is impossible since then the coefficient of $t^4$ vanishes and the coefficient of $t^6$ does not vanish ($(\cala_1)^4=0$ and $(\cala_1)^6\not=0$ is impossible). Thus it remains the cases $s+1=3,N=2$ and $s+1=2,N=3$ for which one has polynomial growth \cite{art-sch:1987}. These latter cases are the object of \cite{art-sch:1987} and we shall describe examples with exponential growth in Section \ref{EX}.

\begin{proposition}\label{M3}
Two 3-regular $(N+1)$-linear forms $w$ and $w'$ on $\mathbb K^{s+1}$ correspond to isomorphic graded algebra $\cala(w,N)$ and $\cala(w',N) $if and only if they belong to the same $GL(s+1,\mathbb K)$-orbits.
\end{proposition}
\noindent \underbar{Proof}. If $w'=w\circ L$ for $L\in GL(s+1,\mathbb K)$, the fact that the corresponding algebras are isomorphic is immediate since then $L$ is just a linear change of generators.\\

Assume now that the graded algebras are isomorphic. Then in degree 1 this isomorphism gives an element $L$ of $GL(s+1,\mathbb K)$ such that in components one has
\begin{equation}
W'_{\alpha_0\alpha_1\dots \alpha_N}=K^\alpha_{\alpha_0} L^{\beta_0}_\alpha L^{\beta_1}_{\alpha_1} \dots L^{\beta_N}_{\alpha_N} W_{\beta_0\beta_1\dots \beta_N}
\end{equation}
for some $K\in GL(s+1,\mathbb K)$ since in view of $(i)$ the relations are linearily independent. Using the property $(ii)$ for $w'$ and for $w$, one gets
\begin{equation}
(Q_{w'})^\alpha_{\alpha_N} W'_{\alpha\alpha_0\dots \alpha_{N-1}}=(K^{-1}L^{-1}Q_wL)^\alpha_{\alpha_N} K^\beta_{\alpha_0}W'_{\alpha\beta\alpha_1\dots \alpha_{N-1}}
\end{equation}
from which it follows by using the property $(iii)$ for $w'$
\begin{equation}
Q_{w'}(L^{-1}Q_wL)^{-1}K=K=\lambda \bbbone
\end{equation}
for some $\lambda \in \mathbb K$. Since $K$ is invertible, one has $Q_{w'}=Q_{w\circ L}$ and $\lambda\not= 0$, and thus $w'=\lambda w\circ L$. i.e. $w'=w\circ L$ by replacing $L$ by $\lambda^{-\frac{1}{N+1}}L$.$\square$\\

\section{Homogeneous algebras associated with multilinear forms}
In this section $m$ and $N$ are two integers such that $m\geq N\geq 2$ and $w$ is a preregular $m$-linear form on $\mathbb K^{s+1}$ with components $W_{\lambda_1\dots \lambda_m}=w(e_{\lambda_1},\dots,e_{\lambda_m})$ in the canonical basis $(e_\lambda)$ of $\mathbb K^{s+1}$. Let $\cala=\cala(w,N)$ be the $N$-homogeneous algebra generated by the $s+1$ elements $x^\lambda$ ($\lambda \in \{0,\dots,s\}$) with the relations
\begin{equation}
W_{\lambda_1\dots \lambda_{m-N}\mu_1\dots \mu_N} x^{\mu_1}\dots x^{\mu_N}=0
\end{equation}
$\lambda_i\in\{0,\dots,s\}$, that is one has $\cala=A(E,R)$ with $E=\oplus_\lambda\mathbb K x^\lambda=\cala_1$ and $R=\sum_{\lambda_i}\mathbb K W_{\lambda_1\dots \lambda_{m-N}\mu_1\dots \mu_N} x^{\mu_1}\otimes \dots \otimes x^{\mu_N}\subset E^{\otimes^N}$. Define $\calw_n\subset E^{\otimes^n}$ for $m\geq n\geq 0$ by 
\begin{equation}
\calw_n=\left\{ \begin{array}{l}
\sum_{\lambda_i} \mathbb K W_{\lambda_1\dots\lambda_{m-n}\mu_1\dots \mu_n} x^{\mu_1}\otimes \dots \otimes x^{\mu_n}\ \ \mbox{if}\ \ m\geq n\geq N\\
\\
E^{\otimes^n}\ \ \ \mbox{if} \ \ \ N-1\geq n\geq 0
\end{array}
\right.
\end{equation}
and consider the sequence
\begin{equation}
\label{Sm}
0\rightarrow \cala \otimes \calw_m \stackrel{d}{\rightarrow} \cala \otimes \calw_{m-1} \stackrel{d}{\rightarrow} \dots \stackrel{d}{\rightarrow} \cala\otimes E \stackrel{d}{\rightarrow} \cala\rightarrow 0
\end{equation}
of free left $\cala$-modules where the homomorphisms $d$ are induced by the homomorphisms of $\cala\otimes E^{\otimes^{n+1}}$ into $\cala\otimes E^{\otimes^n}$ defined by $a\otimes (v_0\otimes v_1\otimes \dots \otimes v_n)\mapsto av_0\otimes (v_1\otimes \dots \otimes v_n)$ for $n\geq 0$, $a\in \cala$ and $v_i\in E=\cala_1$.

\begin{proposition}\label{KSm}
Sequence $(\ref{Sm})$ is a sub-$N$-complex of $K(\cala)$ $($the Koszul $N$-complex of $\cala)$.
\end{proposition}

\noindent \underbar{Proof.} By the property $(ii)$ of $w$ the relations $W_{\lambda_1\dots \lambda_{m-N}\mu_1\dots \mu_N}x^{\mu_1}\dots x^{\mu_N}=0$ are equivalent to $W_{\lambda_{r+1}\dots \lambda_{m-N}\mu_1\dots \mu_N \lambda_1 \dots \lambda_r} x^{\mu_1}\dots x^{\mu_N}=0$ for $m-N\geq r\geq 0$. It follows that $\calw_n\subset E^{\otimes^{n-N-r}}\otimes R\otimes E^{\otimes^r}$ for any $r$ such that $n-N\geq r\geq 0$ so $\calw_n\subset \cap_r E^{\otimes^{n-N-r}}\otimes R\otimes E^{\otimes^r}=\left (\cala^!_n\right)^\ast$ and therefore $\cala\otimes \calw_n\subset K_n(\cala)$ for $n\geq N$. The equalities $\cala\otimes \calw_n=K_n(\cala)$ for $N-1\geq n\geq 0$ are obvious. This implies the result. $\square$\\

In the proof of the above proposition we have shown in particular that one has $\calw_m\subset (\cala^!_m)^\ast$ so that $w\in (\cala^!_m)^\ast$. It follows that $w$ composed with the canonical projection $\cala^!\rightarrow \cala^!_m$ onto degree $m$ defines a linear form $\omega_w$ on $\cala^!$. On the other hand one can write $Q_w\in GL(s+1,\mathbb K)=GL(E^\ast)=GL(\cala^!_1)$. With these notations one has the following proposition.

\begin{proposition} \label{pF}
The element $Q_w$ of $\GL(\cala^!_1) (=\GL (s+1,\mathbb K))$ induces an automorphism $\sigma_w$ of $\cala^!$ and one has $\omega_w(xy)=\omega_w(\sigma_w(y)x)$ for any $x,y \in \cala^!$. Considered as an element $Q^w$ of $GL(\cala_1)=GL(E)$, the transposed $Q^t_w=Q^w$ of $Q_w$ induces an automorphism $\sigma^w$ of $\cala$.
\end{proposition}

\noindent \underbar{Proof}. $Q_w$ induces an automorphism $\tilde \sigma_w$ of degree 0 of the tensor algebra $T(E^\ast)$. Let $\tilde x\in E^{\ast\otimes^N}$ be in $R^\perp$ i.e. $\tilde x=\rho^{\mu_1\dots \mu_N}e_{\mu_1} \otimes \dots \otimes e_{\mu_N}$ with $W_{\lambda_1\dots \lambda_{m-N} \mu_1\dots\mu_N}\rho^{\mu_1\dots \mu_N}=0$, $\forall \lambda_i$. The invariance of $w$ by $Q_w$ implies 
\[
Q^{\rho_1}_{\lambda_1}\dots Q^{\rho_{m-N}}_{\lambda_{m-N}} W_{\rho_1\dots \rho_{m-N} \nu_1\dots \nu_N} Q^{\nu_1}_{\mu_1}\dots Q^{\nu_N}_{\mu_N}\rho^{\mu_1\dots \mu_N}=0
\]
 i.e. $W_{\lambda_1\dots \lambda_{m-N} \nu_1\dots \nu_N} Q^{\nu_1}_{\mu_1} \dots Q^{\nu_N}_{\mu_N} \rho^{\mu_1\dots\mu_N}=0$, $\forall \lambda_i$, which means $\tilde \sigma_w(\tilde x)\in R^\perp$. Thus one has $\tilde\sigma_w(R^\perp)=R^\perp$ which implies that $\tilde \sigma_w$ induces an automorphism $\sigma_w$ of the $N$-homogeneous algebra $\cala^!$. The property $\omega_w(xy)=\omega_w(\sigma_w(y)x)$ for $x,y\in \cala^!$ is then just a rewriting of the property $(ii)$ of $w$ (i.e. the $Q_w$-twisted cyclicity and its consequences given by (3.2)).The last statement of the proposition follows again from the invariance of $w$ by $Q_w$ which implies that one has $(Q^w)^{\otimes^N}(R)\subset R$. $\square$\\
 
 \begin{theorem} \label{Fr}
 The subset $\cali$ of $\cala^!$ defined by 
 \[
 \cali=\{y\in \cala^!\vert \omega_w(xy)=0,\ \ \ \forall x\in \cala^!\}
 \]
 is a graded two-sided ideal of $\cala^!$ and the quotient algebra $\calf(w,N)=\cala^!/\cali$ equipped with the linear form induced by $\omega_w$ is a graded Frobenius algebra.\end{theorem}
 
 \noindent \underbar{Proof}. By its very definition, $\cali$ is a left ideal. It follows from $\omega_w(xy)=\omega_w(\sigma_w(y)x)$ that $\cali$ is also a right ideal, so it is a two-sided ideal. By construction one has $\calf(w,N)=\calf=\oplus^m_{p=0}\calf_p$ with $\dim(\calf_m)=1$ and the pairing induced by $(x,y)\mapsto \omega_w(xy)$ is nondegenerate and is a Frobenius pairing on $\calf$.$\square$\\
 
 One has $\dim(\calf_0)=\dim(\calf_m)=1,\ \dim(\calf_1)=\dim(\calf_{m-1})=s+1$ and of course $\dim(\calf_p)=\dim(\calf_{m-p})$ for $p\in \{0,\dots,m\}$. The automorphism $\sigma_w$ induces an automorphism $\sigma$ of $\calf$ and one has $xy=\sigma(y)x$ for $x\in \calf_p$ and $y\in\calf_{m-p}$, $m\geq p\geq 0$.\\
 
 Notice that if $L\in GL(s+1,\mathbb K)$ then $\cala(w,N)$ and $\cala(w\circ L,N)$ are isomorphic $N$-homogeneous algebras.\\
  
 In the following we let $^w\cala$ denote the $(\cala,\cala)$-bimodule which coincides with $\cala$ as right $\cala$-module and is such that the structure of left $\cala$-module is given by left multiplication by $(-1)^{(m-1)n}(\sigma^w)^{-1}(a)$ for $a\in \cala_n$. One has the following result in the quadratic case $\cala=\cala(w,2)$.
 
 \begin{proposition}\label{Vol}
 Assume that $N=2$, that is that $\cala$ is the quadratic algebra $\cala=\cala(w,2)$. Then $\bbbone \otimes w$ is canonically a nontrivial $^w\cala$-valued Hochschild $m$-cycle on $\cala$, that is one has $\bbbone \otimes w \in Z_m(\cala,^w\cala)$ with $\bbbone \otimes w \not\in B_m(\cala,^w\cala)$.
 \end{proposition}
 
 \noindent \underbar{Proof}. The $m$-linear form $w$ on $\mathbb K^{s+1}$ identifies canonically with an element of $E^{\otimes^m}=\cala_1^{\otimes^m}
\subset \cala^{\otimes^m}$, i.e. one can write $w\in \cala^{\otimes^m}$. By 
interpreting $\bbbone \in \cala$ as an element of $^w\cala$ one can consider that $\bbbone\otimes w$ is an $^w\cala$-valued Hochschild $m$-chain. The Hochschild boundary of $\bbbone\otimes w$ reads
\[
\begin{array}{lll}
b(\bbbone\otimes w) & = & W_{\lambda_1\dots \lambda_m}x^{\lambda_1}\otimes \dots \otimes^{\lambda_m}\\
&  +  & \sum^{m-1}_{k=1} (-1)^k \bbbone \otimes W_{\lambda_1\dots \lambda_m} x^{\lambda_1} \otimes \dots \otimes x^{\lambda_k}x^{\lambda_k+1}\otimes\dots x^{\lambda_{m}}\\
& - & (Q^{-1}_w)^{\lambda_m}_\lambda W_{\lambda_1\dots \lambda_m} x^\lambda \otimes x^{\lambda_1}\otimes\dots \otimes x^{\lambda_{m-1}}
\end{array}
\]
The sum of the first and of the last term vanishes by $Q_w$-cyclicity while each of the other terms vanishes since the relations $W_{\lambda_1\dots \lambda_{m-2}\mu\nu} x^\mu x^\nu=0$ are equivalent to $W_{\lambda_1\dots \lambda_r \mu\nu \lambda_{r+1}\dots \lambda_{m-2}}x^\mu x^\nu=0$ again by $Q_w$-cyclicity. Therefore one has $b(\bbbone\otimes w)=0$. By using the fact that the Hochschild boundary preserves the total $\cala$-degree it is easy to see that $\bbbone\otimes w$ cannot be a boundary.$\square$\\
Thus in the quadratic case $N=2$, if $Q_w=(-1)^{m-1}$ then $\bbbone\otimes w$ represents the analog of a differential $m$-form i.e. an element of $H_m(\cala,\cala)=HH_m(\cala)$; if $Q_w$ is different of $(-1)^{m-1}$, this is a twisted version of a differential $m$-form. We shall come back later to this interpretation in the Koszul-Gorenstein case where $\bbbone\otimes w$ plays the role of a volume element.

\begin{theorem}\label{KGD}
Let $\cala$ be a $N$-homogeneous algebra generated by $s+1$ elements $x^\lambda$ $(\lambda \in \{0,\dots,s\})$ which is Koszul of finite global dimension $D$ and which is Gorenstein. Then $\cala=\cala(w,N)$ for some preregular $m$-linear form $w$ on  $\mathbb K^{s+1}$ which is such that if $N\geq 3$ then one has $m=Np+1$ and $D=2p+1$ for some integer $p\geq 1$ while for $N=2$ one has $m=D$. Under these conditions $w$ is unique up to a multiplicative factor in $\mathbb K\backslash \{0\}$.
\end{theorem}

\noindent \underbar {Proof}. The Koszul resolution of the trivial left $\cala$-module $\mathbb K$ ends as ( 
\cite{ber-mdv-wam:2003}, see also Appendix 1)
\begin{equation}
\dots \stackrel{d}{\rightarrow} \cala \otimes \cala^{!\ast}_N \stackrel{d^{N-1}}{\rightarrow} \cala\otimes \cala_1 \stackrel{d}{\rightarrow} \cala \rightarrow \mathbb K \rightarrow 0
\end{equation}
so the Gorenstein property implies that it starts as
\begin{equation}
0\rightarrow \cala \otimes \cala^{!\ast}_m \stackrel{d}{\rightarrow} \cala \otimes \cala^{!\ast}_{m-1} \stackrel{d^{N-1}}{\rightarrow} \dots
\end{equation}
with 
\begin{equation}
\dim (\cala^{!\ast}_m)=1
\end{equation}
and
\begin{equation}
\dim (\cala^{!\ast}_{m-1})=\dim(\cala_1)=s+1
\end{equation}
for some $m\geq N$ which corresponds to the $D$-th term. This implies that $m=D$ for $N=2$ and that for $N\geq 3$ one has $m=Np+1$ and $D=2p+1$ for some integer $p\geq 1$. Let $w$ be a generator of the 1-dimensional subspace $\cala^{!\ast}_m$ of $\cala^{\otimes^m}_1$. Since $\cala_1$ identifies canonically with the dual vector space of $\mathbb K^{s+1}$, $w$ is (canonically) a $m$-linear form on $\mathbb K^{s+1}$. For $0\leq k\leq m$ and $\theta \in \cala^!_k$, one defines $\theta w\in \cala^{!\ast}_{m-k}$ by setting
\begin{equation}
(\theta w)(\alpha)=\langle w,\alpha\theta\rangle
\end{equation}
for $\alpha\in \cala^!_{m-k}$. The mapping $\theta\mapsto \theta w$ defines a left $\cala^!$-module homomorphism $\bar \Phi$ of $\cala^!$ into $\cala^{!\ast}$ with $\bar\Phi(\cala^!_k)\subset \cala^{!\ast}_{m-k}$ for $k\in\{0,\dots,m\}$ and the Gorenstein property implies that $\bar\Phi$ induces the linear isomorphisms
\begin{equation}\label{Iso}
\bar\Phi:\cala^!_{\nu_N(p)}\simeq \cala^{!\ast}_{\nu_N(D-p)}
\end{equation}
for $p\in \{0,\dots,D\}$ where $\nu_N(2\ell)=N\ell$ and $\nu_N(2\ell+1)=N\ell+1$ for $\ell\in \mathbb N$, \cite{ber-mar:2006} (Theorem 5.4 of \cite{ber-mar:2006}). The isomorphisms (\ref{Iso}) for $p=1$ and $p=D-1$ imply that $w$ is a preregular $m$-linear form on $\mathbb K^{s+1}$  while the isomorphism (\ref{Iso}) for $p=D-2$ implies that the relations of $\cala$ read
\[
W_{\lambda_1\dots \lambda_{m-N}\mu_1\dots \mu_N} x^{\mu_1}\dots x^{\mu_N}=0
\]
with $W_{\lambda_1\dots \lambda_m}=w(e_{\lambda_1},\dots,e_{\lambda_m})=\langle w,e_{\lambda_1}\dots e_{\lambda_m}\rangle$ and hence that one has $\cala=\cala(w,N)$. $\square$\\

Notice that under the assumptions of Theorem \ref{KGD}, the Koszul resolution of the trivial left $\cala$-module $\mathbb K$ reads
\[
0\rightarrow \cala\otimes \calw_m\stackrel{d}{\rightarrow} \cala\otimes \calw_{m-1} \stackrel{d^{N-1}}{\rightarrow} \dots \stackrel{d}{\rightarrow} \cala\otimes \calw_N \stackrel{d^{N-1}}{\rightarrow} \cala\otimes E \stackrel{d}{\rightarrow} \cala \rightarrow \mathbb K \rightarrow 0
\]
i.e.
\begin{equation}
0\rightarrow \cala \otimes \calw_{\nu_N(D)} \stackrel{d'}{\rightarrow} \dots \stackrel{d'}{\rightarrow} \cala \otimes \calw_{\nu_N(k)} \stackrel{d'}{\rightarrow} \cala \otimes \calw_{\nu_N(k-1)} \stackrel{d'}{\rightarrow} \dots \stackrel{d'}{\rightarrow} \cala\rightarrow \mathbb K \rightarrow 0
\end{equation}
where $d'=d^{N-1}:\cala\otimes \calw_{\nu_N(2\ell)}\rightarrow \cala \otimes \calw_{\nu_N(2\ell-1)}$ and $d'=d:\cala\otimes \calw_{\nu_N(2\ell+1)}\rightarrow \cala\otimes \calw_{\nu_N(2\ell)}$ and that one has
\begin{equation}
\dim(\calw_{\nu_N(k)})=\dim(\calw_{\nu_N(D-k)})
\end{equation}
for any $0\leq k\leq D$. Thus $\cala\otimes \calw_m=\cala\otimes \calw_{\nu_N(D)}=\cala\otimes w$ so one sees that $\bbbone \otimes w$ is the generator of the top module of the Koszul resolution which also leads to an interpretation of $\bbbone \otimes w$ as volume form.

\section{Examples}\label{EX}

\subsection{Yang-Mills algebra} \label{subsec1}

Let $g_{\mu\nu}$ be the components of a symmetric nondegenerate bilinear form on $\mathbb K^{s+1}$. The Yang-Mills algebra \cite{ac-mdv:2002b}, \cite{nek:2003} is the cubic algebra $\cala$ generated by the $s+1$ elements $x^\lambda$ ($\lambda \in \{0,\dots s\}$) with the $s+1$ relations
\begin{equation}\label{YM}
g_{\lambda\mu}[x^\lambda,[x^\mu,x^\nu]]=0
\end{equation}
for $\nu\in \{0,\dots, s\}$.\\

It was claimed in \cite{ac-mdv:2002b} that this algebra is Koszul of global dimension 3 and is Gorenstein. The relations (\ref{YM}) can be rewritten as
\begin{equation}
(g_{\rho\lambda} g_{\mu\nu}+g_{\rho\nu}g_{\lambda\mu}-2g_{\rho\mu}g_{\lambda\nu})x^\lambda x^\mu x^\nu=0
\end{equation}
(for $\rho\in \{0,\dots, s\}$) and one verifies that the 4-linear form $w$ on $\mathbb K^{s+1}$ with components
\begin{equation}
W_{\rho\lambda \mu \nu}=g_{\rho\lambda} g_{\mu\nu}+g_{\rho\nu} g_{\lambda\mu}-2g_{\rho\mu}g_{\lambda\nu}
\end{equation}
is 3-regular with $Q_w=\bbbone$. So it is invariant by cyclic permutations and one has $\cala=\cala(w,3)$ with the notations of the previous sections. It is easy to see that $w$ does not only satisfy the condition $(iii)$ but satisfies the stronger condition $(iii)'$. This implies (and is equivalent to the fact) that the dual cubic algebra $\cala^!$ is Frobenius.

\subsection{Super Yang-Mills algebra}

With the same conventions as in \ref{subsec1}, the super Yang-Mills algebra 
\cite{ac-mdv:2007} is the cubic algebra $\tilde \cala$ generated by the $s+1$ elements $x^\lambda$ with the $s+1$ relations

\begin{equation}\label{SYM}
g_{\lambda\mu}[x^\lambda,\{x^\mu,x^\nu\}]=0
\end{equation}
for $\nu\in \{0,\dots, s\}$ and where $\{A,B\}=AB+BA$. As pointed out in 
\cite{ac-mdv:2007} the relations (\ref{SYM}) can be equivalently written as
\begin{equation}
[g_{\lambda\mu}x^\lambda x^\mu, x^\nu]=0
\end{equation}
which means that the quadratic element $g_{\lambda\mu}x^\lambda x^\mu$ is central.\\
It was claimed in \cite{ac-mdv:2007} that this algebra is Koszul of global dimension 3 and is Gorenstein.\\
The relations (\ref{SYM}) can  be rewritten as
\begin{equation}
(g_{\rho\lambda} g_{\mu\nu} -g_{\rho\nu} g_{\lambda\mu})x^\lambda x^\mu x^\nu
\end{equation}
and one verifies that the 4-linear form $\tilde w$ on $\mathbb K^{s+1}$ with components
\begin{equation}
\tilde W_{\rho\lambda\mu\nu}=g_{\rho\lambda}g_{\mu\nu} - g_{\rho\nu} g_{\lambda\mu}
\end{equation}
is 3-regular with $Q_{\tilde w}=-\bbbone$ and satisfies the stronger condition $(iii)'$. Thus one has $\tilde \cala=\cala(\tilde w,3)$ and the dual cubic algebra $\tilde\cala^!$ is Frobenius.\\

The Poincar\'e series of the Yang-Mills algebra and of the super Yang-Mills algebra coincide and are given by 
\begin{equation}
\frac{1}{1-(s+1)t+(s+1)t^3-t^4}=\left(\frac{1}{1-t^2}\right)\left(\frac{1}{1-(s+1)t+t^2}\right)
\end{equation}
which is a particular case of the formula (\ref{S3}). It follows that these algebras have exponential growth for $s\geq 2$. For $s+1=2$ these are particular cubic Artin-Schelter algebras \cite{art-sch:1987} (see also \cite{mdv-pop:2002}).\\

An elegant powerful proof of the Koszulity of the Yang-Mills and the super Yang-Mills algebras is given in \cite{kri-vdb:2010}.

\subsection{The algebras $\cala(\varepsilon, N)$ for $s+1\geq N \geq 2$}

Assume that $s+1\geq N\geq 2$ and let $\varepsilon$ be the completely antisymmetric $(s+1)$-linear form on $\mathbb K^{s+1}$ such that $\varepsilon(e_0,\dots,e_s)=1$, where $(e_\lambda)$ is the canonical basis of $\mathbb K^{s+1}$. It is clear that $\varepsilon$ is preregular with $Q_\varepsilon=(-1)^s\bbbone$ and that furthermore it satisfies $(iii)$ whenever $s+1\geq 3$.\\
The $N$-homogeneous algebra $\cala(\varepsilon, N)$ generated by the $s+1$ elements $x^\lambda$ with the relations
\begin{equation}
\varepsilon_{\alpha_0\dots \alpha_{s-N}\lambda_1\dots \lambda_N} x^{\lambda_1}\dots x^{\lambda_N}
\end{equation}
(for $\alpha_i\in \{0,\dots,s\}$) was introduced in \cite{ber:2001a} where it was shown that it is Koszul of finite global dimension. It was then shown in \cite{ber-mar:2006}  that $\cala(\varepsilon,N)$ is Gorenstein if and only if either $N=2$ or $N>2$ and $s=Nq$ for some integer $q\geq 1$. The case $N=2$ corresponds to the polynomial algebra with $s+1$ indeterminates which is of global dimension $s+1$ while in the case $N>2$ and $s=Nq$ the $N$-homogeneous algebra $\cala(\varepsilon, N)$ which is then Koszul and Gorenstein has global dimension $D=2q+1$.\\
In the general case if $N>2$ the dual $N$-homogeneous algebra $\cala(\varepsilon, N)^!$ cannot be Frobenius since the ideal $I_\varepsilon$ always contains the non trivial quadratic elements
\begin{equation}
e_\lambda e_\mu +e_\mu e_\lambda
\end{equation}
and is in fact generated by these elements in the Koszul-Gorenstein case i.e. when $s=Nq$ with $q\geq 1$. In this latter case $\cala(\varepsilon,N)^!/I_\varepsilon$ is the exterior algebra $\wedge\mathbb K^{s+1}$ over $\mathbb K^{s+1}$ which is the dual $\cala(\varepsilon,2)^!$ ($=\wedge\mathbb K^{s+1}$) of the quadratic algebra $\cala(\varepsilon,2)$ generated by the $x^\lambda$ ($\lambda\in \{0,\dots,s\}$) with the relations
\begin{equation}
x^\lambda x^\mu=x^\mu x^\lambda
\end{equation}
(which coincides with the polynomial algebra with $s+1$ indeterminates). \\
One thus recovers by this process the quadratic relations implying the original $N$-homogeneous ones with $N>2$ (for $s=Nq$ with $q\geq 1$). One may wonder whether there is a lesson to extract from this example : Namely starting from a Koszul-Gorenstein $N$-homogeneous algebra $\cala$ with $N$-homogeneous dual $\cala^!$ which is not Frobenius, is there a $N_0$-homogeneous $\cala_0$ with $\cala^!_0$ Frobenius such that the relation of $\cala$ are implied by the relations of $\cala_0$ ? Notice that the Koszul-Gorenstein algebras $\cala(\varepsilon,N)$ (for $s=Nq,N\geq 3, q\geq 1$) have exponential growth.

\subsection{The algebra $\cala_\ug$}

Let us now discuss in the present context the algebra $\cala_\ug$ introduced in \cite{ac-mdv:2002a} and analyzed in details in 
 \cite{ac-mdv:2003}
 and \cite{ac-mdv:2008}. The algebra $\cala_\ug$, which corresponds to a noncommutative 4-plane, is the quadratic algebra generated by the 4 generators $x^\lambda$ ($\lambda\in \{0,1,2,3\}$) with the relations
 \begin{equation}
 \cos (\varphi_0-\varphi_k)[x^0,x^k]=i\sin (\varphi_\ell - \varphi_m) \{x^\ell,x^m\}
 \end{equation}
 \begin{equation}
 \cos(\varphi_\ell-\varphi_m)[x^\ell,x^m]=i\sin (\varphi_0-\varphi_k)\{x^0,x^k\}
 \end{equation}
 for any cyclic permutation $(k,\ell,m)$ of (1,2,3) and where $\{A,B\}=AB+BA$ as before. The parameter $\ug$ being the element
 \begin{equation}
 \ug = (e^{i(\varphi_1-\varphi_0)},e^{i(\varphi_2-\varphi_0)}, e^{i(\varphi_3-\varphi_0)})
 \end{equation}
 of $T^3$. This algebra is Koszul of global dimension 4 and is Gorenstein (so an example with $N=2$ and $s+1=D=4$) whenever none of these six relations becomes trivial and one has the non trivial Hochschild cycle \cite{ac-mdv:2002a}
 \begin{eqnarray}
 w_\ug = \tilde{\mbox{ch}}_{3/2}(U_\ug) & =  &-\sum_{\alpha, \beta, \gamma, \delta} \varepsilon_{\alpha\beta\gamma\delta} \cos (\varphi_\alpha-\varphi_\beta + \varphi_\gamma-\varphi_\delta) x^\alpha \otimes x^\beta \otimes x^\gamma \otimes x^\delta\nonumber \\
& + & i\sum_{\mu,\nu} \sin (2(\varphi_\mu-\varphi_\nu)) x^\mu \otimes x^\nu \otimes x^\mu \otimes x^\nu
\end{eqnarray}
which may be considered as a 4-linear form on $\mathbb K^4$ which is preregular. Notice that the components $W_{\rho\lambda\mu\nu}$ of $w_\ug=W_{\rho\lambda\mu\nu} x^\rho\otimes x^\lambda\otimes x^\mu\otimes x^\nu$ can be written as
\begin{equation}
W_{\rho\lambda\mu\nu}=-\cos(\varphi_\rho-\varphi_\lambda + \varphi_\mu-\varphi_\nu)\varepsilon_{\rho\lambda\mu\nu}+i\sin (\varphi_\rho-\varphi_\lambda+\varphi_\mu -\varphi_\nu)\delta_{\rho\mu}\delta_{\lambda\nu}
\end{equation}
for $ \rho, \lambda, \mu, \mu \in \{0, 1, 2, 3\}$.\\

One can check that one has $\cala_\ug =\cala(w_\ug,2)$ and, as explained in \cite{ac-mdv:2002a}, one has $Q_{w_\ug}=-\bbbone$ and $\bbbone \otimes w_\ug$ is a Hochschild 4-cycle ($\in Z_4(\cala_\ug,\cala_\ug)$) which is a particular case of the corresponding more general result of Section 5 in the quadratic case. \\
As pointed out in \cite{ac-mdv:2002a}, for generic values of the parameter $\ug$ the algebra $\cala_\ug$ is isomorphic to the Sklyanin algebra \cite{skl:1982} which has been studied in detail in \cite{smi-sta:1992} from the point of view of regularity. 

\section {Semi-cross product (twist)}

In this section we investigate semi-cross products of algebras of type $\cala(w,N)$ by homogeneous automorphisms of degree 0 and we describe the corresponding transformations of the $Q_w$. We first recall the definition and the properties of the semi-cross product using the notations of \cite{ac-mdv:2008}. This notion has been introduced and analyzed in \cite{art-tat-vdb:1991} where it is refered to as twisting and used to reduce and complete the classification of regular algebras of dimension 3 of \cite{art-sch:1987}. In \cite{ac-mdv:2008} this notion was used to reduce the moduli space of noncommutative 3-spheres.\\

Let $\cala=\otimes_{n\in \mathbb N}\cala_n$ be a graded algebra and let $\alpha$ be an automorphism of $\cala$ which is homogeneous of degree 0. The {\sl semi-cross product $\cala(\alpha)$ of $\cala$ by $\alpha$} is the graded vector space $\cala$ equipped with the bilinear product $\bullet_\alpha=\bullet$ defined by 
\[
a\bullet b = a \alpha^n(b)
\]
for $a\in \cala_n$ and $b\in \cala$. This product is associative and satisfies $\cala_m\bullet \cala_n\subset \cala_{m+n}$ so $\cala(\alpha)$ is a graded algebra which is unital whenever $\cala$ is unital. The following result is extracted from \cite{art-tat-vdb:1991}.

\begin{proposition}
Let $\cala$, $\alpha$ and $\cala(\alpha)$ be as above.\\
$(i)$ The global dimensions of $\cala$ and $\cala(\alpha)$ coincide.\\
$(ii)$ Let $\beta$ be an automorphism of $\cala$ which is homogeneous of degree 0 and which commutes with $\alpha$. Then $\beta$ is also canonically an automorphism of $\cala(\alpha)$ and one has
\[
\cala(\alpha)(\beta)=\cala(\alpha\circ \beta)
\]
which implies in particular that $\cala(\alpha)(\alpha^{-1})=\cala$.
\end{proposition}
In fact the category of graded right $\cala$-modules and the category of graded right $\cala(\alpha)$-modules are canonically isomorphic which implies $(i)$. For $(ii)$ we refer to \cite{art-tat-vdb:1991} (see also in 
\cite{ac-mdv:2008}). \\
If $\cala$ is $N$-homogeneous then $\cala(\alpha)$ is also $N$-homogeneous and if $R\subset \cala^{\otimes^N}_1$ denotes the space of relations of $\cala$, then the space of relations of $\cala(\alpha)$ is given by \cite{ac-mdv:2008}
\begin{equation}
R(\alpha)=(Id\otimes \alpha^{-1}\otimes \dots \otimes \alpha^{-(N-1)})R
\label{A1}
\end{equation}
with obvious notations. Concerning the stability of the Koszul property and of the Gorenstein property, the following result was proved in \cite{pot:2006}.
\begin{proposition}
Let $\cala$ be a $N$-homogeneous algebra and $\alpha$ be an homogeneous automorphism of degree 0 of $\cala$.\\
$(i)$ $\cala(\alpha)$ is Koszul if and only if $\cala$ is Koszul.\\
$(ii)$ $\cala(\alpha)$ is Koszul of (finite) global dimension $D$ and Gorenstein if and only if $\cala$ is Koszul of global dimension $D$ and Gorenstein.
\end{proposition}

If $\cala$ is a $N$-homogeneous algebra, an automorphism $\alpha$ of degree 0 of $\cala$ is completely specified by its restriction $\alpha\restriction \cala_1$ to $\cala_1$.\\

Let $w$ be a preregular $m$-linear form on $\mathbb K^{s+1}$ with $m\geq N\geq 2$ and let us consider the $N$-homogeneous algebra $\cala=\cala(w,N)$. We denote by $GL_w$ the subgroup of $GL(s+1,\mathbb K)$ of the elements $L\in GL(s+1,\mathbb K)$ which preserve $w$, i.e. such that
\begin{equation}
w\circ L =w
\label{A2}
\end{equation}
It is clear that each $L\in GL_w$ determines an automorphism $\alpha^{(L)}$ of degree 0 of $\cala$ for which $\alpha^{(L)}\restriction \cala_1$ is the transposed $L^t$ of $L$. Furthermore, it follows from (\ref{A1}) and (\ref{A2}) that the semi-cross product of $\cala$ by $\alpha^{(L)}$ is given by
\begin{equation}
\cala(\alpha^{(L)})=\cala(w^{(L)},N)
\label{A3}
\end{equation}
where the components $W^{(L)}_{\lambda_1\dots \lambda_m}$ of the $m$-linear form $w^{(L)}$ are given by
\begin{equation}
W^{(L)}_{\lambda_1\dots\lambda_m}=W_{\lambda_1\lambda'_2\dots\lambda'_m}(L^{-1})^{\lambda'_2}_{\lambda_2} (L^{-2})^{\lambda'_3}_{\lambda_3}\dots
(L^{-(m-1)})^{\lambda'_m}_{\lambda_m}
\label{A4}
\end{equation}
in terms of the components $W_{\mu_1\dots\mu_m}$ of $w$. The $m$-linear form $w^{(L)}$ is again preregular with $Q_{w^{(L)}}$ given by
\begin{equation}
Q_{w^{(L)}}=L^{-1}Q_w L^{-(m-1)}
\label{A5}
\end{equation}
as verified by using (\ref{A4}) and (\ref{A2}).\\

The fact that $\cala(w,N)$ and $\cala(w\circ L,N)$ are isomorphic $N$-homogeneous algebras for $L\in GL(s+1,\mathbb K)$ implies in view of Theorem 11 that, in the study of Koszul-Gorenstein algebras of finite global dimension, one can simplify $Q_w$ by using $Q_w\mapsto L^{-1}Q_wL=Q_{w\circ L}$ ($L\in GL(s+1,\mathbb K)$; e.g. one can assume that $Q_w$ is in Jordan normal form. On the other hand, since the construction of the semi-cross product is very explicit and since it preserves the global dimension (Proposition 12) and the Koszul-Gorenstein property (Proposition 13), it is natural to simplify further  $Q_w$ via $Q_w\mapsto L^{-1}Q_wL^{-(m-1)}$ with $L\in GL_w$ (formula (A5)). In many cases one can find a $m$-th root of $\pm Q_w$ in $GL_w$, that is an element $L\in GL(s+1,\mathbb K)$ such that $w\circ L=w,\ [Q_w,L]=0$ and $L^m=\pm Q_w$. In such a case one can restrict attention to $Q_w=\pm \bbbone$, i.e. to $w$ which is $\pm$ cyclic, by semi-cross product (twist). There are however some cases where this is not possible (in fact there are cases where $GL_w$ is a small discrete group).\\

It is worth noticing here that a more general notion of twisted graded algebras which is quite optimal has been introduced and analyzed in \cite{zha:1996}. Several results given here are valid for this more general twisting as pointed out there.

\section{Actions of quantum groups}

Although it is evident that in this section we only need 1-site  nondegenerate multilinear forms (see Section 3), we shall stay in the context of Section 5. That is we let $m$ and $N$ be two integers with $m\geq N\geq 2$ and $w$ be a pre-regular $m$-linear form on $\mathbb K^{s+1}$ with components $W_{\lambda_1\dots \lambda_m}=w(e_{\lambda_1},\dots,e_{\lambda_m})$ in the canonical basis $(e_\lambda)$ of $\mathbb K^{s+1}$. As pointed out in Section 2 and in more details in Appendix 2, in the case $m=N=2$, that is when $w$ is a nondegenerate bilinear form $b$, there is a Hopf algebra $\calh(b)$ and a natural coaction of $\calh(b)$ on $\cala(b,2)$ or, in dual terms there is a quantum group acting on the quantum space corresponding to $\cala(b,2)$. Our aim in this section is to generalize this and to define a Hopf algebra $\calh$ (in fact several generically) which coacts on $\cala(w,N)$ for general $m\geq N\geq 2$. For the cases where the $\cala(w,N)$ are Artin-Schelter regular algebras, there are the closely related works  \cite{ewe-ogi:1994} and \cite{pop:2006}. Here however we merely concentrate on the ``$SL$-like" aspect (instead of the ``$GL$-like" one). This is also closely related to the quantum $SU$ of \cite{wor:1988} and the quantum $SL$ of \cite{bic:2001}.\\

By the 1-site nondegeneracy property of $w$, there is (at least one) a $m$-linear form $\tilde w$ on the dual of $\mathbb K^{s+1}$ with components $\tilde W^{\lambda_1\dots \lambda_m}$ in the dual basis of $(e_\lambda)$ such that one has
\begin{equation}
\tilde W^{\alpha\gamma_1\dots\gamma_{m-1}}W_{\gamma_1\dots \gamma_{m-1}\beta}= \delta^\alpha_\beta\label{7.1}
\end{equation}
for $\alpha, \beta\in \{0,\dots,s\}$. Let $\calh(w,\tilde w)$ be the unital associative algebra generated by the $(s+1)^2$ elements $u^\alpha_\beta$ ($\alpha,\beta\in \{0,\dots,s\}$) with the relations
\begin{equation}
W_{\alpha_1\dots\alpha_m} u^{\alpha_1}_{\beta_1}\dots u^{\alpha_m}_{\beta_m}=W_{\beta_1\dots \beta_m}\bbbone\label{7.2}
\end{equation}
and
\begin{equation}
\tilde W^{\beta_1\dots\beta_m}u^{\alpha_1}_{\beta_1}\dots u^{\alpha_m}_{\beta_m}=\tilde W^{\alpha_1\dots\alpha_m}\bbbone\label{7.3}
\end{equation}
where $\bbbone$ is the unit of $\calh(w,\tilde w)$. One has the following result.

\begin{theorem}
There is a unique structure of Hopf algebra on $\calh(w,\tilde w)$ with coproduct $\Delta$, counit $\varepsilon$ and antipode $S$ such that

\begin{eqnarray}
\Delta(u^\mu_\nu) & = & u^\mu_\lambda \otimes u^\lambda_\nu\label{7.4}\\
\varepsilon(u^\mu_\nu) & = & \delta^\mu_\nu\label{7.5}\\
S(u^\mu_\nu) & = & \tilde W^{\mu\lambda_1\dots \lambda_{m-1}}u^{\rho_1}_{\lambda_1}\dots u^{\rho_{m-1}}_{\lambda_{m-1}} W_{\rho_1\dots \rho_{m-1}\nu}\label{7.6}
\end{eqnarray}
for $\mu,\nu\in \{0,\dots,s\}$. The product and the unit being the ones of $\calh(w,\tilde w$).
\end{theorem} 

\noindent \underbar{Proof}. The structure  of bialgebra with (\ref{7.4}) and (\ref{7.5}) is more or less classical. The fact that $S$ defines an antipode follows from $S(u^\mu_\lambda)u^\lambda_\nu=\delta^\mu_\nu$ and $u^\mu_\lambda S(u^\lambda_\nu)=\delta^\mu_\nu$ which are immediate consequences of (\ref{7.1}), (\ref{7.2}), (\ref{7.3}) and (\ref{7.6}).$\square$

\begin{proposition}
There is a unique algebra-homomorphism 
\[
\Delta_L:\cala(w,N)\rightarrow \calh(w,\tilde w)\otimes \cala(w,N)
\]
 such that
\begin{equation}
\Delta_L(x^\mu)=u^\mu_\nu \otimes x^\nu\label{7.7}
\end{equation}
for $\mu\in \{0,\dots,s\}$. This endows $\cala(w,N)$ of a structure of $\calh(w,\tilde w)$-comodule.
\end{proposition}

\noindent \underbar{Proof}. One has $W_{\alpha_1\dots \alpha_{m-N} \mu_1\dots\mu_N}(u^{\alpha_1}_{\lambda_1}\otimes \bbbone)\dots (u^{\alpha_{m-N}}_{\lambda_{m-N}}\otimes \bbbone)\Delta_L x^{\mu_1}\dots \Delta_L x^{\mu_N}= \bbbone \otimes W_{\lambda_1\dots \lambda_{m-N}\mu_1\dots \mu_N} x^{\mu_1} \dots x^{\mu_N}=0$. By contracting on the right-hand side with $(S(u^{\lambda_{m-N}}_{\nu_{m-N}})\otimes \bbbone)\dots (S(u^{\lambda_1}_{\nu_1})\otimes \bbbone)$, this is equivalent to 
\[
W_{\nu_1\dots \nu_{m-N}\mu_1\dots \mu_N}\Delta_L x^{\mu_1}\dots \Delta_Lx^{\mu_N}=0
\]
for $\nu_i\in \{0,\dots,s\}$. The fact that $\Delta_L$ induces a structure of $\calh(w,\tilde w)$-comodule is straightforward.$\square$\\

It is worth noticing that the preregularity of $w$ implies (in view of the $Q_w$-cyclicity) that one has
\begin{equation}
\tilde W^{\lambda\gamma_2\dots \gamma_m}W_{\mu\gamma_2\dots \gamma_m}=(Q^{-1}_w)^\lambda_\mu
\end{equation}
which generalizes a formula valid for $m=N=2$.\\

The dual object of the Hopf algebra $\calh(w,\tilde w)$ is a quantum group which acts on the quantum space corresponding (dual object) to the algebra $\cala(w,N)$.\\

The above theorem and the above proposition correspond to the theorem and the proposition of Appendix 2 for the case $m=N=2$. There is however a notable difference in the cases $m\geq N\geq 2$ which is that for $m=N=2$ that is when $w=b$, a nondegenerate bilinear form, then $\tilde w$ is unique under condition (\ref{7.1}) and coincides with $b^{-1}$ (see Appendix 2), $\tilde w=b^{-1}$. In the general case, given $w$ there are several $\tilde w$ satisfying (\ref{7.1}) and thus several Hopf algebras $\calh(w,\tilde w)$. Some choice of $\tilde w$ can be better than other in the sense that $\calh(w,\tilde w)$ can be bigger or can have bigger commutative quotient. For instance, in the case of Example 6.3 of Section 6 where $w=\varepsilon$, the natural choice for $\tilde w$ is $\tilde \varepsilon$ with components $\frac{(-1)^s}{s!}\varepsilon^{\lambda_0\dots \lambda_s}$ where $\varepsilon^{\lambda_0\dots \lambda_s}$ is completely antisymmetric with $\varepsilon^{0\ 1\dots s}=1$; in this case, $\calh(\varepsilon, \tilde \varepsilon)$ is commutative \cite{wor:1988} and coincides with the Hopf algebra of representative functions on $SL(s+1,\mathbb K)$.\\

In fact the right quantum group relevant here corresponds to the universal Hopf algebra $\calh(w)$ preserving $w$ which has been recently studied in \cite{bic-mdv:2013} where an explicit finite presentation by generators and relations is given for it. It turns out that for $w=\varepsilon$, $\calh(\varepsilon)$ is not commutative whenever $m\geq 3$ and therefore that the canonical Hopf algebra map $\calh(\varepsilon)\rightarrow \calh(\varepsilon, \tilde\varepsilon)$ is not injective for $m\geq 3$, (see in \cite{bic-mdv:2013}).

\section{Further prospect}

As pointed out in the introduction it was already shown in \cite{bon-pol:1994} that the quadratic algebras which are Koszul of finite global dimension $D$ and which are Gorenstein are determined by multilinear forms ($D$-linear forms). Furthermore the connection with a generalization of volume forms is also apparent in \cite{bon-pol:1994}. This corresponds to the case $N=2$ of Theorem 11. The argument of \cite{bon-pol:1994} is that the Koszul dual algebra $\cala^!$ of a quadratic algebra $\cala$ which is Koszul of finite global dimension and Gorenstein is a graded Frobenius algebra which is generated in degree 1 and that such an algebra is completely characterized by a multilinear form on the ($D$-dimensional) space of generators $\cala^!_1$ of $\cala^!$. Here, the argument is slightly different and works in the other way round (in the quadratic case). Indeed, Theorem 9 combined with Theorem 11 imply that in the quadratic case the Koszul dual $\cala^!$ of a Koszul Gorenstein algebra $\cala$ of finite global dimension is Frobenius. This points however to the interesting observation that for a $N$-homogeneous algebra $\cala$ which is Koszul of finite global dimension and Gorenstein there are two graded Frobenius algebras which can be extracted from $\cala^!$. These two graded Frobenius algebras coincide with $\cala^!$ in the quadratic case but are distinct whenever $N\geq 3$. The first one is the Yoneda algebra $E(\cala)=\Ext^\bullet_\cala(\mathbb K, \mathbb K)$ which can be obtained by truncation from $\cala^!$ as explained in \cite{ber-mar:2006} while the second one is the quotient $\cala^!/\cali$ of Theorem 9. The Yoneda algebra $E(\cala)$ has been considered by many authors as the generalization of the Koszul dual for nonquadratic graded algebras. However the quotient $\cala^!/\cali$ is also of great interest for homogeneous algebras and deserves further attention (see e.g. the discussion at the end of \S 6.3) .\\

The converse of Theorem 5 which was stated as a conjecture in the last version of this paper (v3) is unfortunately wrong. One has the following counterexample which can be found in \cite{art-sch:1987}. Let $\cala$ be the algebra generated by 3 elements $x,y,z$ with relations $x^2+yz=0,\> y^2+zx=0,\> xy=0$. Then $\cala=\cala(w,2)$ with $w$ given by
\[
w=x^{\otimes^3}+y^{\otimes^3}+x\otimes y\otimes z + y\otimes z\otimes x + z\otimes x\otimes y
\]
which is a 3-regular 3-linear form on $\mathbb K^3$. It turns out that $\cala$ is nevertheless not Koszul, see e.g. in \cite{mdv:2010}. However the following results are worth noticing.

\begin{proposition}\label{CK3}
Let $w$ be a preregular $(N+1)$-linear form on $\mathbb K^{s+1}$ and $\cala=\cala(w,N)$. Then the following conditions are equivalent.\\
$(i)$ $w$ is 3-regular\\
$(ii)$ The dual vector space $(\cala^!_{N+1})^\ast$ of $\cala^!_{N+1}$ is given by  $(\cala^!_{N+1})^\ast=\mathbb Kw$\\
$(iii)$ The Koszul complex $\calk(\cala,\mathbb K)$ of $\cala$ is the sequence 
\[
0\rightarrow \cala\otimes w \stackrel{d}{\longrightarrow} \cala\otimes R \stackrel{d^{N-1}}{\longrightarrow} \cala\otimes E \stackrel{d}{\longrightarrow} \cala \rightarrow 0
\]
$(iiii)$ The Koszul $N$-complex $K(\cala)$ of $\cala$ is the sequence
\[
0\rightarrow \cala\otimes w \stackrel{d}{\longrightarrow} \cala\otimes R \stackrel{d}{\longrightarrow} \cala\otimes E^{\otimes^{N-1}} \stackrel{d}{\longrightarrow} \dots \stackrel{d}{\longrightarrow} \cala\otimes E \stackrel{d}{\longrightarrow}\cala\rightarrow 0
\]
where $E=\oplus_\lambda \mathbb K x^\lambda=\cala_1$, $R=\oplus_\mu\mathbb K W_{\mu\mu_1\dots \mu_N}x^{\mu_1}\otimes \dots \otimes x^{\mu_N}\subset E^{\otimes^N}$ and where, in $(i)$, $(ii)$ and $(iii)$, $w$ is identified with $W_{\mu_0\dots \mu_N} x^{\mu_0}\otimes \dots \otimes x^{\mu_N}\in E^{\otimes^{N+1}}$.
\end{proposition}
\noindent \underbar{Proof}. 
\begin{itemize}
\item
$(iii)\Leftrightarrow (iiii)$. This follows from the definitions and from the fact that 
\end{itemize}
$\cala^!_{N+2}=0$ implies $\cala^!_n=0$ for $n\geq N+2$ in view of the associativity of the product of $\cala^!$.
\begin{itemize}
\item
$(i)\Leftrightarrow (ii)$. Let $a\in (\cala^!_{N+1})^\ast=(R\otimes E) \cap (E\otimes R)$ then
\end{itemize} 
\[
a=W_{\lambda_0\dots \lambda_{N-1}\rho} L^\rho_{\lambda_N} x^{\lambda_0}\otimes \dots \otimes x^{\lambda_N}=M^\sigma_{\lambda_0}W_{\sigma\lambda_1\dots\lambda_N}x^{\lambda_0}\otimes \dots \otimes x^{\lambda_N}
\]
so one has
\begin{equation}
W_{\lambda_0\dots \lambda_{N-1}\rho} L^\rho_{\lambda_N}=M^\sigma_{\lambda_0}W_{\sigma\lambda_1\dots \lambda_N}, \> \> \> \forall \lambda_i
\label{a}
\end{equation}
and conversely any solution of (\ref{a}) defines an element $a$ of $(\cala^!_{N+1})^\ast$.
By the preregularity property (twisted cyclicity) of $w$, (\ref{a})  is equivalent to 
\begin{equation}
(Q^{-1}_w)_{\lambda_0}^\alpha L^\beta_\alpha (Q_w)^\tau_\beta W_{\tau\lambda_1\dots \lambda_N}=M^\sigma_{\lambda_1} W_{\lambda_0\sigma \lambda_2\dots \lambda_N},\>\>\> \forall \lambda_i
\label{b} 
\end{equation}
and $a=kw$ $(k\in \mathbb K)$ is then equivalent (1-site nondegeneracy) to $L=M=k\bbbone$ (or equivalently $Q_wLQ_w^{-1}=M=k\bbbone$). Since $a\in (R\otimes E)\cap (E\otimes R)$ is arbitrary this implies $(i)\Leftrightarrow (ii)$.
\begin{itemize}
\item
$(iiii)\Rightarrow (ii)$. This follows from $K_{N+1}(\cala)=\cala\otimes (\cala^!_{N+1})^\ast$.
\item
$(i)\Rightarrow (iii)$. In order to complete the proof it is sufficient to show that
\end{itemize}
if $w$ is 3-regular then $(\cala^!_{N+2})^\ast=0$. So assume now that $w$ is 3-regular and let $a\in (\cala^!_{N+2})^\ast=(E^{\otimes^2}\otimes R)\cap (E\otimes R\otimes E)\cup (R\otimes E^{\otimes^2})$. One has $a=A^\lambda_{\lambda_0\lambda_1}W_{\lambda\lambda_2\dots \lambda_{N+1}}x^{\lambda_0}\otimes \dots \otimes x^{\lambda_{N+1}}$ with
\begin{equation}
A^\lambda_{\lambda_0 \lambda_1} W_{\lambda\lambda_2\dots \lambda_{N+1}}=B^\rho_{\lambda_0\lambda_{N+1}}W_{\rho\lambda_1\dots \lambda_N}
\label{c}
\end{equation}
\begin{equation}
B^\rho_{\lambda_0\lambda_{N+1}}W_{\rho\lambda_1\dots \lambda_N}=C^\sigma_{\lambda_N\lambda_{N+1}}W_{\sigma\lambda_0\dots \lambda_{N-1}}
\label{d}
\end{equation}
for any $\lambda_i$. By the 3-regularity of $w$, equation (\ref{c}) implies 
\[
A^\lambda_{\nu\mu}=\left(Q^{-1}_w\right)^\tau_\mu B^\lambda_{\nu\tau}=K_\nu \delta^\lambda_\mu
\]
while equation (\ref{d}) implies
\[
B^\lambda_{\nu\mu} = \left(Q^{-1}_w\right)^\tau_\nu  C^\lambda_{\tau\mu}=L_\mu \delta^\lambda_\nu
\]
and therefore one has
\begin{equation}
K_{\lambda_0}W_{\lambda_1\dots \lambda_{N+1}}=W_{\lambda_0\dots \lambda_N} L_{\lambda_{N+1}}
\label{(e)}
\end{equation}
for any $\lambda_i$. Since $w$ is 1-site nondegenerate this implies $K=L=0$ so $A=B=C=0$ and therefore $a=0$. Thus if $w$ is 3-regular then $(\cala^!_{N+2})^\ast=0$.$\square$

\begin{corollary}
Let $w$ be a 3-regular $(N+1)$-linear form on $\mathbb K^{s+1}$ and assume that the $N$-homogeneous algebra $\cala=\cala(w,N)$ is Koszul. Then $\cala$ is Koszul of global dimension 3 and is Gorenstein.
\end{corollary}

\noindent \underbar{Proof}. From the last proposition it follows that one has the (Koszul) minimal projective resolution 
\[
0\rightarrow \cala\otimes w \stackrel{d}{\rightarrow} \cala\otimes R \stackrel{d^{N-1}}{\longrightarrow} \cala\otimes E \stackrel{d}{\rightarrow} \cala \rightarrow \mathbb K \rightarrow 0
\]
of the trivial left $\cala$-module $\mathbb K$. The Gorenstein property is then equivalent to the twisted cyclicity of $w$ (property $(ii)$ of Definition 2); this is the same argument as the one used in \cite{art-sch:1987}. Another way to prove this result is to use the corollary 5.12 of \cite{ber-mar:2006} since it is clear that the 1-site nondegenerate property of $w$ implies here that the Yoneda algebra $E(\cala)$ is Frobenius.~$\square$\\

In fact it is the same here to assume that $\cala(w,N)$ is of global dimension 3 as to assume that it is Koszul and one has the following result.\\

\begin{theorem}
Let $\cala$ be a connected graded algebra which is finitely generated in degree 1 and finitely presented with relations of degree $\geq 2$. Then $\cala$ has global dimension $D=3$ and is Gorenstein if and only if it is Koszul of the form $\cala=\cala(w,N)$ for some 3-regular $(N+1)$-linear form $w$ on $\mathbb K^{s+1}$ $(s+1=\dim \cala_1)$.
\end{theorem}

It is possible to give higher dimensional generalization of the 3-regularity, namely $D$-regularity for (preregular) $D$-linear forms ($D\geq N=2$) and\linebreak[4] $(2q+1)$-regularity for (preregular) $(Nq+1)$-linear forms ($N\geq 2$). However the cases $D=4$ ($N=2$) and $D=5$ are already very cumbersome.

\subsubsection*{Acknowledgements}
It is a pleasure to thank Roland Berger and Alain Connes for their kind interest and for their suggestions as well as, for this version, Andrea Solotar, Paul Smith and Michel Van den Bergh for their critical comments and advices.

\newpage

\section*{Appendix 1: Homogeneous algebras}
\setcounter{section}{10}
\setcounter{equation}{0}
A {\sl homogeneous algebra of degree $N$ or $N$-homogeneous algebra} is an algebra of the form 
\[
\cala = A(E,R)=T(E)/(R)
\]
where $E$ is a finite-dimensional vector space, $R$ is a linear subspace of $E^{\otimes^N}$ and where $(R)$ denotes the two-sided ideal of the tensor algebra $T(E)$ of $E$ generated by $R$. The algebra $\cala$ is naturally a connected graded algebra with graduation induced by the one of $T(E)$. To $\cala$ is associated another $N$-homogeneous algebra, {\sl its dual}  $\cala^!=A(E^\ast, R^\perp)$ with $E^\ast$ denoting the dual vector space of $E$ and $R^\perp\subset E^{\otimes^N\ast}=E^{\ast \otimes^N}$ being the annihilator of $R$, \cite{ber-mdv-wam:2003}. The $N$-complex $K(\cala)$ of left $\cala$-modules is then defined to be
\begin{equation}
\dots \stackrel{d}{\rightarrow} \cala\otimes \cala^{!\ast}_{n+1} \stackrel{d}{\rightarrow} \cala\otimes \cala^{!\ast}_{n}\stackrel{d}{\rightarrow} \dots \stackrel{d}{\rightarrow}\cala \rightarrow 0
\label{eq7.1}
\end{equation}
where $\cala^{!\ast}_n$ is the dual vector space of the finite-dimensional vector space $\cala^!_n$ of the elements of degree $n$ of $\cala^!$ and where $d:\cala\otimes \cala^{!\ast}_{n+1}\rightarrow \cala\otimes \cala^{!\ast}_n$ is induced by the map $a\otimes (e_1\otimes \dots \otimes e_{n+1})\mapsto ae_1 \otimes (e_2\otimes \dots \otimes e_{n+1})$ of $\cala\otimes E^{\otimes^{n+1}}$ into $\cala\otimes E^{\otimes^n}$, remembering that $\cala^{!\ast}_n\subset E^{\otimes^n}$, (see \cite{ber-mdv-wam:2003}). This $N$-complex will be refered to as the {\sl Koszul $N$-complex of} $\cala$.  In (\ref{eq7.1}) the factors $\cala$ are considered as left $\cala$-modules. By considering $\cala$ as right $\cala$-module and by exchanging the factors one obtains the $N$-complex $\tilde K(\cala)$ of right $\cala$-modules
\begin{equation}
\dots \stackrel{\tilde d}{\rightarrow} \cala^{!\ast}_{n+1} \otimes \cala \stackrel{\tilde d}{\rightarrow} \cala^{!\ast}_n\otimes \cala \stackrel{\tilde d}{\rightarrow} \dots \stackrel{\tilde d}{\rightarrow} \cala \rightarrow 0
\label{eq7.2}
\end{equation}
where now $\tilde d$ is induced by $(e_1\otimes \dots \otimes e_{n+1})\otimes a \mapsto (e_1\otimes \dots \otimes e_n)\otimes e_{n+1}a$. Finally one defines two $N$-differentials $d_{\Lg}$ and $d_{\Rg}$ on the sequence of $(\cala,\cala)$-bimodules, i.e. of left $\cala\otimes \cala^{opp}$-modules, $(\cala\otimes \cala^{!\ast}_n \otimes \cala)_{n\geq 0}$ by setting $d_\Lg =d\otimes I_\cala$ and $d_\Rg = I_\cala\otimes \tilde d$ where $I_\cala$ is the identity mapping of $\cala$ onto itself. For each of these $N$-differentials $d_\Lg$ and $d_\Rg$ the sequences
\begin{equation}
\dots \stackrel{d_\Lg,d_\Rg}{\rightarrow} \cala \otimes \cala^{!\ast}_{n+1}\otimes \cala \stackrel{d_\Lg,d_\Rg}{\rightarrow} \cala\otimes \cala^{!\ast}_n \otimes \cala \stackrel{d_\Lg, d_\Rg}{\rightarrow}\dots
\label{eq7.3}
\end{equation}
are $N$-complexes of left $\cala\otimes \cala^{opp}$-modules and one has
\begin{equation}
d_\Lg d_\Rg = d_\Rg d_\Lg
\label{eq7.4}
\end{equation}
which implies that
\begin{equation}
d^N_\Lg -d^N_\Rg = (d_\Lg -d_\Rg)\left ( \sum^{N-1}_{p=0} d^p_\Lg d^{N-p-1}_\Rg\right) = \left( \sum^{N-1}_{p=0} d^p_\Lg d^{N-p-1}_\Rg\right) (d_\Lg -d_\Rg) =0
\label{eq7.5}
\end{equation}
in view of $d^N_\Lg=d^N_\Rg=0$.\\

As for any $N$-complex \cite{mdv:1998a} one obtains from $K(\cala)$ ordinary complexes $C_{p,r}(K(\cala))$, {\sl the contractions of} $K(\cala)$, by putting together alternatively $p$ and $N-p$ arrows $d$ of $K(\cala)$. Explicitely $C_{p,r}(K(\cala))$ is given by 
\begin{equation}
\dots \stackrel{d^{N-p}}{\rightarrow} \cala \otimes \cala^{!\ast}_{Nk+r}\stackrel{d^p}{\rightarrow} \cala\otimes \cala^{!\ast}_{Nk-p+r}\stackrel{d^{N-p}}{\rightarrow}\cala\otimes \cala^{!\ast}_{N(k-1)+r}\stackrel{d^p}{\rightarrow}\dots
\label{eq7.6}
\end{equation}
for $0\leq r< p\leq N-1$, \cite{ber-mdv-wam:2003} . These are here chain complexes of free left $\cala$-modules. As shown in \cite{ber-mdv-wam:2003} the complex $C_{N-1,0}(K(\cala))$ coincides with the {\sl Koszul complex} of \cite{ber:2001a}; this complex will be denoted by $\calk(\cala,\mathbb K)$ in the sequel. That is one has
\begin{equation}
\calk_{2m}(\cala,\mathbb K)=\cala\otimes\cala^{!\ast}_{Nm},\ \ \ \calk_{2m+1}(\cala,\mathbb K)=\cala\otimes \cala^{!\ast}_{Nm+1}
\label{eq7.7}
\end{equation}
for $m\geq 0$, and the differential is $d^{N-1}$ on $\calk_{2m}(\cala,\mathbb K)$ and $d$ on $\calk_{2m+1}(\cala,\mathbb K)$. If $\calk(\cala, \mathbb K)$ is acyclic in positive degrees then $\cala$ will be said to be a {\sl Koszul algebra}. It was shown in \cite{ber:2001a} and this was confirmed by the analysis of \cite{ber-mdv-wam:2003} that this is the right generalization for $N$-homogeneous algebra of the usual notion of Koszulity for quadratic algebras  \cite{man:1987}, \cite{lod:1999}. One always has $H_0(\calk(\cala,\mathbb K))\simeq \mathbb K$ and therefore if $\cala$ is Koszul, then one has a free resolution $\calk(\cala,\mathbb K)\rightarrow \mathbb K \rightarrow 0$ of the trivial left $\cala$-module $\mathbb K$, that is the exact sequence
\begin{equation}
\dots \stackrel{d^{N-1}}{\rightarrow} \cala \otimes \cala^{!\ast}_{N+1}\stackrel{d}{\rightarrow}\cala\otimes R\stackrel{d^{N-1}}{\rightarrow}\cala\otimes E\stackrel{d}{\rightarrow} \cala \stackrel{\varepsilon}{\rightarrow} \mathbb K \rightarrow 0
\label{eq7.8}
\end{equation}
of left $\cala$-modules where $\varepsilon$ is the projection on degree zero. This resolution is a minimal projective resolution of the $\cala$-module $\mathbb K$  in the graded category  \cite{ber:2005} which will be refered to as the {\sl Koszul resolution of $\mathbb K$}.\\

One defines now the chain complex of free $\cala\otimes\cala^{opp}$-modules $\calk(\cala,\cala)$ by setting
\begin{equation}
\calk_{2m}(\cala,\cala)=\cala\otimes \cala^{!\ast}_{Nm}\otimes \cala,\ \ \ \calk_{2m+1}(\cala,\cala)=\cala\otimes \cala^{!\ast}_{Nm+1}\otimes \cala
\label{eq7.9}
\end{equation}
for $m\in \mathbb N$ with differential $\delta'$ defined by
\begin{equation}
\delta'=d_\Lg -d_\Rg :\calk_{2m+1}(\cala,\cala)\rightarrow \calk_{2m}(\cala,\cala)
\label{eq7.10}
\end{equation}
\begin{equation}
\delta'=\sum^{N-1}_{p=0} d^p_\Lg d^{N-p-1}_\Rg:\calk_{2(m+1)}(\cala,\cala)\rightarrow \calk_{2m+1}(\cala,\cala)
\label{eq7.11}
\end{equation}
the property $\delta^{\prime 2}=0$ following from (\ref{eq7.5}). {\sl This complex is acyclic in positive degrees if and only if $\cala$ is Koszul}, that is if and only if $\calk(\cala,\mathbb K)$ is acyclic in positive degrees,  \cite{ber:2001a} and \cite{ber-mdv-wam:2003} . One always has the obvious exact sequence
\begin{equation}
\cala\otimes E\otimes \cala \stackrel{\delta'}{\rightarrow} \cala\otimes \cala \stackrel{\mu}{\rightarrow}\cala\rightarrow 0
\label{eq7.12}
\end{equation}
of left $\cala\otimes \cala^{opp}$-modules where $\mu$ denotes the product of $\cala$. It follows that if $\cala$ is a Koszul algebra then $\calk(\cala,\cala)\stackrel{\mu}{\rightarrow}\cala\rightarrow 0$ is a free resolution of the $\cala\otimes \cala^{opp}$-module $\cala$ which will be refered to as {\sl the Koszul resolution of} $\cala$. This is a minimal projective resolution of the $\cala\otimes \cala^{opp}$-module $\cala$  in the graded category \cite{ber:2005}.\\

Let $\cala$ be a Koszul algebra and let $\calm$ be a $(\cala,\cala)$-bimodule considered as a right $\cala\otimes \cala^{opp}$-module. Then, by interpreting the $\calm$-valued Hochschild homology $H(\cala,\calm)$ as $H_n(\cala,\calm)=\Tor_n^{\cala\otimes \cala^{opp}}(\calm,\cala)$  {\cite{car-eil:1973},  the complex $\calm \otimes_{\cala\otimes \cala^{opp}}\calk(\cala,\cala)$ computes the $\calm$-valued Hochschild homology of $\cala$, (i.e. its homology is the ordinary $\calm$-valued Hochschild homology of $\cala$). We shall refer to this complex as {\sl the small Hochschild complex of} $\cala$ with coefficients in $\calm$ and denote it by $\cals(\cala, \calm)$. It reads
\begin{equation}
\dots \stackrel{\delta}{\rightarrow}\calm\otimes \cala^{!\ast}_{N(m+1)}\stackrel{\delta}{\rightarrow} \calm \otimes \cala^{!\ast}_{Nm+1}\stackrel{\delta}{\rightarrow}\calm\otimes \cala^{!\ast}_{Nm} \stackrel{\delta}{\rightarrow}\dots
\label{eq7.13}
\end{equation}
where $\delta$ is obtained from $\delta'$ by applying the factors $d_L$ to the right of $\calm$ and the factors $d_R$ to the left of $\calm$.\\

Assume that $\cala$ is a Koszul algebra of finite global dimension $D$. Then the Koszul resolution of $\mathbb K$ has length $D$, i.e. $D$ is the largest integer such that $\calk_D(\cala,\mathbb K)\not=0$. By construction, $D$ is also the greatest integer such that $\calk_D(\cala,\cala)\not= 0$ so the free $\cala\otimes \cala^{opp}$-module resolution of $\cala$ has also length $D$. Thus one verifies in this case the general statement of  \cite{ber:2005} namely that the global dimension is equal to the Hochschild dimension. Applying then the functor $\Hom_\cala(\bullet, \cala)$ to $\calk(\cala,\mathbb K)$ one obtains the cochain complex $\call(\cala,\mathbb K)$ of free right $\cala$-modules
\[
0\rightarrow \call^0(\cala,\mathbb K)\rightarrow \dots \rightarrow \call^D(\cala,\mathbb K) \rightarrow 0
\]
where $\call^n(\cala,\mathbb K)=\Hom_\cala(\calk_n(\cala,\mathbb K),\cala)$. The Koszul algebra $\cala$ is {\sl Gorenstein} iff $H^n(\call(\cala,\mathbb K))=0$ for $n<D$ and $H^D(\call(\cala,\mathbb K))=\mathbb K$ (= the trivial right $\cala$-module). This is clearly a generalisation of the classical Poincar\'e duality and this implies a precise form of Poincar\'e duality between Hochschild homology and Hochschild cohomology \cite{ber-mar:2006}, \cite{vdb:1998}, \cite{vdb:2002}. In the case of the Yang-Mills algebra and its deformations which are Koszul Gorenstein cubic algebras of global dimension 3, this Poincar\'e duality gives isomorphisms
\begin{equation}
H_k(\cala,\calm)= H^{3-k}(\cala,\calm),\>\>\> k\in \{0,1,2,3\}
\label{Pd}
\end{equation}
between the Hochschild homology and the Hochschild cohomology with coefficients in a bimodule $\calm$. This follows from the fact that in these cases one has $Q_w=\bbbone$.

\section*{Appendix 2: The quantum group of a nondegenerate bilinear form}

Let $b$ be a non degenerate bilinear form on $\mathbb K^{s+1}$ with components $B_{\mu\nu}=b(e_\mu,e_\nu)$ in the canonical basis $(e_\lambda)_{\lambda\in\{0,\dots,s\}}$. The matrix elements $B^{\mu\nu}$ of the inverse $B^{-1}$ of the matrix $B=(B_{\mu\nu})$ of components of $b$ are the components of a nondegenerate bilinear form $b^{-1}$ on the dual vector space of $\mathbb K^{s+1}$ in the dual basis of $(e_\lambda)$. Let $\calh(b)$ be the unital associative algebra generated by the $(s+1)^2$ elements $u^\mu_\nu$ ($\mu,\nu\in \{0,\dots,s\}$) with the relations
\begin{equation}
B_{\alpha\beta} u^\alpha_\mu u^\beta_\nu=B_{\mu\nu}\bbbone\ \ \ (\mu,\nu\in \{0,\dots,s\})
\end{equation}
and
\begin{equation}
B^{\mu\nu} u^\alpha_\mu u^\beta_\nu = B^{\alpha\beta}\bbbone\ \ \ (\alpha,\beta \in \{0,\dots,s\})
\end{equation}
where $\bbbone$ denotes the unit of $\calh(b)$. One has the following \cite{mdv-lau:1990}.
\begin{theorem}
There is a unique structure of Hopf algebra on $\calh(b)$ with coproduct $\Delta$, counit $\varepsilon$ and antipode $S$ such that
\begin{eqnarray}
\Delta (u^\mu_\nu) & = & u^\mu_\lambda \otimes u^\lambda_\nu\\
\varepsilon (u^\mu_\nu) & = &  \delta^\mu_\nu\\
S(u^\mu_\nu) & = & B^{\mu\alpha}B_{\beta\nu} u^\beta_\alpha
\end{eqnarray}
for $\mu,\nu\in \{0,\dots,s\}$. 
The product and the unit being the ones of $\calh(b)$. 
\end{theorem}
The proof is straightforward and the dual object of the Hopf algebra $\calh(b)$ is called {\sl the quantum group of the nondegenerate bilinear form $b$} ;  $\calh(b)$ corresponds to the Hopf algebra of ``representative functions" on this quantum group.
\begin{proposition}
Let $\cala=\cala(b,2)$ be the $($quadratic$)$ algebra of Section 2 and $\calh=\calh(b)$ be the above Hopf algebra. There is a unique algebra-homomorphism $\Delta_L:\cala\rightarrow \calh\otimes \cala$ such that
\begin{equation}
\Delta_L(x^\lambda)=u^\lambda_\mu \otimes x^\mu
\end{equation} 
for$\lambda\in\{0,\dots,s\}$. This endows $\cala$ of a structure of $\calh$-comodule.
\end{proposition}
Thus the quantum group of $b$ ``acts" on the quantum space corresponding to $\cala$.\\

 Let $q\in \mathbb K$ with $q\not=0$ be such that 
\begin{equation}
B^{\alpha\beta} B_{\alpha\beta} + q + q^{-1}=0
\label{eq8.7}
\end{equation}
then the linear endomorphisms $R_\pm$ of $\mathbb K^{s+1}\otimes K^{s+1}$ defined by
\begin{equation}
(R_+)^{\alpha\beta}_{\mu\nu}=\delta^\alpha_\mu \delta^\beta_\nu+qB^{\alpha\beta}B_{\mu\nu},\ \ (R_-)^{\alpha\beta}_{\mu\nu}=\delta^\alpha_\mu \delta^\beta_\nu+q^{-1} B^{\alpha\beta} B_{\mu\nu}
\end{equation}
satisfy the Yang-Baxter relation
\begin{equation}
(R_\pm \otimes \bbbone)(\bbbone \otimes R_\pm)(R_\pm\otimes \bbbone)=(\bbbone \otimes R_\pm) (R_\pm\otimes \bbbone)(\bbbone \otimes R_\pm): (\mathbb K^{s+1})^{\otimes^3}\rightarrow (\mathbb K^{s+1})^{\otimes^3}
\end{equation} and $(R_+-1)(R_++q^2)=0$, $(R_--1)(R_-+q^{-2})=0$.\\

Let $\varepsilon_q$ ($q\not=0$) be the nondegenerate bilinear form on $\mathbb K^2$ with matrix of components $\left ( \begin{array}{cc}
0 & -1\\
q & 0
\end{array} \right)$. Then $\cala(\varepsilon_q,2)=\cala_q$ corresponds to the Manin plane (see in Section 2) whereas $\calh(\varepsilon_q)=\calh_q$ corresponds to the quantum group $SL_q(2,\mathbb K)$.\\

 One has the following result of  \cite{bic:2003b} concerning the representations of the quantum group of the nondegenerate bilinear form $b$ on $\mathbb K^{s+1}$.
\begin{theorem}
Let $b$ be a nondegenerate bilinear form on $\mathbb K^{s+1}$ and let $q\in \mathbb K\backslash \{0\}$ be defined by $(\ref{eq8.7})$. Then the category of comodules on $\calh(b)$ is equivalent to the category of comodules on $\calh(\varepsilon_q)=\calh_q$.
\end{theorem}
In other words, in the dual picture, the category of representations of the quantum group of the nondegenerate bilinear form $b$ is equivalent to the category of representations of the quantum group $SL_q(2,\mathbb K)$ with $q$ given by (\ref{eq8.7}).

%

\begin{thebibliography}{10}

\bibitem{art-sch:1987}
M.~Artin and W.F. Schelter.
\newblock Graded algebras of global dimension 3.
\newblock {\em Adv. Math.}, 66:171--216, 1987.

\bibitem{art-tat-vdb:1990}
M.~Artin, J.~Tate, and M.~Van~den Bergh.
\newblock {S}ome algebras associated to automorphisms of elliptic curves. {T}he
  {G}rothendieck {F}estschrift. vol.{ I}.
\newblock {\em Prog. Math.}, 86:33--85, 1990.

\bibitem{art-tat-vdb:1991}
M.~Artin, J.~Tate, and M.~Van~den Bergh.
\newblock Modules over regular algebras of dimension 3.
\newblock {\em Invent. Math.}, 106:335--388, 1991.

\bibitem{ber:2001a}
R.~Berger.
\newblock Koszulity for nonquadratic algebras.
\newblock {\em J. Algebra}, 239:705--734, 2001.

\bibitem{ber:2005}
R.~Berger.
\newblock Dimension de {H}ochschild des alg{\`e}bres gradu{\'e}es.
\newblock {\em C.R. Acad.Sci. Paris, Ser. I}, 341:597--600, 2005.

\bibitem{ber:2009}
R.~Berger.
\newblock Gerasimov's theorem and ${N}$-{K}oszul algebras.
\newblock {\em J. London Soc.}, 79:631--648, 2009.

\bibitem{ber-mdv-wam:2003}
R.~Berger, M.~Dubois-Violette, and M.~Wambst.
\newblock Homogeneous algebras.
\newblock {\em J. Algebra}, 261:172--185, 2003.

\bibitem{ber-mar:2006}
R.~Berger and N.~Marconnet.
\newblock Koszul and {G}orenstein properties for homogeneous algebras.
\newblock {\em Algebras and Representation Theory}, 9:67--97, 2006.

\bibitem{bic:2001}
J.~Bichon.
\newblock Cosovereign {H}opf algebras.
\newblock {\em J. Pure Appl. Algebra}, 157:121--133, 2001.

\bibitem{bic:2003b}
J.~Bichon.
\newblock The representation category of the quantum group of a non-degenerate
  bilinear form.
\newblock {\em Comm. Algebra}, 31:4831--4851, 2003.

\bibitem{bic-mdv:2013}
J.~Bichon and M.~Dubois-Violette.
\newblock The quantum group of a preregular multilinear form.
\newblock {\em Lett. Math. Phys.}, 103:455--468, 2013.

\bibitem{bon-pol:1994}
A.I. Bondal and A.E. Polishchuk.
\newblock Homological properties of associative algebras: The method of
  helices.
\newblock {\em Russian Acad. Sci. Izv. Math.}, 42:219--260, 1994.

\bibitem{car-eil:1973}
H.~Cartan and S.~Eilenberg.
\newblock {\em Homological algebra}.
\newblock Princeton University Press, 1973.

\bibitem{ac:1986a}
A.~Connes.
\newblock Non-commutative differential geometry.
\newblock {\em Publ. IHES}, 62:257--360, 1986.

\bibitem{ac:1994}
A.~Connes.
\newblock {\em Non-commutative geometry}.
\newblock Academic Press, 1994.

\bibitem{ac-mdv:2002a}
A.~Connes and M.~Dubois-Violette.
\newblock Noncommutative finite-dimensional manifolds. {I}. {S}pherical
  manifolds and related examples.
\newblock {\em Commun. Math. Phys.}, 230:539--579, 2002.

\bibitem{ac-mdv:2002b}
A.~Connes and M.~Dubois-Violette.
\newblock Yang-{M}ills algebra.
\newblock {\em Lett. Math. Phys.}, 61:149--158, 2002.

\bibitem{ac-mdv:2003}
A.~Connes and M.~Dubois-Violette.
\newblock Moduli space and structure of noncommutative 3-spheres.
\newblock {\em Lett. Math. Phys.}, 66:91--121, 2003.

\bibitem{ac-mdv:2007}
A.~Connes and M.~Dubois-Violette.
\newblock Yang-{M}ills and some related algebras.
\newblock In {\em Rigorous {Q}uantum {F}ield {Theory}}, volume 251 of {\em
  Progr. Math.}, pages 65--78. Birkha{\"u}ser, 2007.

\bibitem{ac-mdv:2008}
A.~Connes and M.~Dubois-Violette.
\newblock Noncommutative finite-dimensional manifolds. {II}. {M}oduli space and
  structure of noncommutative 3-spheres.
\newblock {\em Commun. {M}ath. {P}hys.}, 281(1):23--127, 2008.

\bibitem{mdv:1998a}
M.~Dubois-Violette.
\newblock $d^{N}=0$: Generalized homology.
\newblock {\em $K$-{T}heory}, 14:371--404, 1998.

\bibitem{mdv:2005}
M.~Dubois-Violette.
\newblock Graded algebras and multilinear forms.
\newblock {\em C. R. Acad. Sci. Paris, Ser. {I}}, 341:719--724, 2005 (see
  ArXiv:math/0509689 v3 for a recent actualization).

\bibitem{mdv:2010}
M.~Dubois-Violette.
\newblock Noncommutative coordinate algebras.
\newblock In E.~Blanchard, editor, {\em Quanta of {M}aths, d{\'e}di{\'e} {\`a}
  {A}. {C}onnes}, volume~14 of {\em Clay Mathematics Proceedings}, pages
  171--199. Clay Mathematics Institute, 2010.

\bibitem{mdv-lau:1990}
M.~Dubois-Violette and G.~Launer.
\newblock The quantum group of a non-degenerate bilinear form.
\newblock {\em Phys. Lett.}, 245B:175--177, 1990.

\bibitem{mdv-pop:2002}
M.~Dubois-Violette and T.~Popov.
\newblock Homogeneous algebras, statistics and combinatorics.
\newblock {\em Lett. Math. Phys.}, 61:159--170, 2002.

\bibitem{ewe-ogi:1994}
H.~Ewen and O.~Ogievetsky.
\newblock Classification of the {GL}(3) quantum matrix groups.
\newblock q-alg/9412009.

\bibitem{irv:1979}
R.S. Irving.
\newblock Prime ideals of {O}re extensions.
\newblock {\em J. Algebra}, 58:399--423, 1979.

\bibitem{kri-vdb:2010}
B.~Kriegk and M.~Van~den Bergh.
\newblock Representations of non-commutative quantum groups.
\newblock ArXiv : 1004.4210.

\bibitem{lod:1999}
J.L. Loday.
\newblock Notes on {K}oszul duality for associative algebras.
\newblock 1999.

\bibitem{lu-pal-wu-zha:2007}
D.M. Lu, J.H. Palmieri, Q.S. Wu, and J.J. Zhang.
\newblock Regular algebras of dimension 4 and their {A}$_\infty$-ext-algebras.
\newblock {\em Duke Math. J.}, 137:537--584, 2007.

\bibitem{man:1987}
Yu.~I. Manin.
\newblock Some remarks on {K}oszul algebras and quantum groups.
\newblock {\em Ann. Inst. Fourier, Grenoble}, 37:191--205, 1987.

\bibitem{man:1988}
Yu.~I. Manin.
\newblock {\em Quantum groups and non-commutative geometry}.
\newblock CRM Universit{\'e} de Montr{\'e}al, 1988.

\bibitem{nek:2003}
N.~Nekrasov.
\newblock Lectures on open strings and noncommutative gauge theories.
\newblock In {\em Les Houches 2001}, NATO Adv. Study Inst., pages 477--495. EDP
  Science, 2003.

\bibitem{pop:2006}
T.~Popov.
\newblock Automorphisms of regular algebras.
\newblock In V.~K. Dobrev, editor, {\em Lie {T}heory and its {A}pplications in
  {P}hysics VI, {V}arna 2006}, 2006.

\bibitem{pot:2006}
A.~Pottier.
\newblock Stabilit\'e de la propri\'et\'e de {K}oszul pour les alg\`ebres
  homog\`enes vis-\`a-vis du produit semi-crois\'e.
\newblock {\em C.R. Acad. Sci. Paris, S{\'e}rie {I}}, 343:161--164, 2006.

\bibitem{pri:1970}
S.B. Priddy.
\newblock Koszul resolutions.
\newblock {\em Trans. Amer. Math. Soc.}, 152:39--60, 1970.

\bibitem{skl:1982}
E.K. Sklyanin.
\newblock Some algebraic structures connected with the {Y}ang-{B}axter
  equation.
\newblock {\em Func. Anal. Appl.}, 16:263--270, 1982.

\bibitem{smi-sta:1992}
S.P. Smith and J.T. Stafford.
\newblock Regularity of the four dimensional {S}klyanin algebra.
\newblock {\em Compos. Math.}, 83:259--289, 1992.

\bibitem{staf:2002}
J.T. Stafford.
\newblock Noncommutative projective geometry.
\newblock {\em ICM}, III:1--3, 2002.

\bibitem{vdb:1998}
M.~Van~den Bergh.
\newblock A relation between {H}ochschild homology and cohomology for
  {G}orenstein rings.
\newblock {\em Proc. Amer. Math. Soc.}, 126:1345--1348, 1998.

\bibitem{vdb:2002}
M.~Van~den Bergh.
\newblock Erratum.
\newblock {\em Proc. Amer. Math. Soc.}, 130:2809--2810, 2002.

\bibitem{wor:1988}
S.L. Woronowicz.
\newblock Tannaka-{K}rein duality for compact matrix pseudogroups. {T}wisted
  {SU}({N}) groups.
\newblock {\em Invent. Math.}, 93:35--76, 1988.

\bibitem{zha:1996}
J.J. Zhang.
\newblock Twisted graded algebras and equivalences of graded categories.
\newblock {\em Proc. London Math. Soc.}, 72:281--311, 1996.

\bibitem{zha:1998}
J.J. Zhang.
\newblock Non-{N}oetherian regular rings of dimension 2.
\newblock {\em Proc. Amer. Math. Soc.}, 126:1645--1653, 1998.

\end{thebibliography}

\end{document}